\documentclass[11pt]{article}
\usepackage{amsmath}
\usepackage{graphicx}
\usepackage{epsfig}
\usepackage{amsfonts}
\usepackage{amssymb}
\usepackage{placeins}
\usepackage{cases}
\usepackage[margin=1in]{geometry}
\usepackage{authblk}
\usepackage[titletoc,toc,title]{appendix}
\usepackage{epstopdf}
\usepackage{changes}
\usepackage{xcolor}
\usepackage{bm}
\usepackage{enumerate}

\def\horizontaldistance{\kern2pt}
\def\verticaldistance{\kern 5pt}

\usepackage{accents}

\title{Identification and Computation of Slow Manifolds Using the Isostable Coordinate System}

\begin{document}
\author[1]{Dan Wilson \thanks{corresponding author:~dwilso81@utk.edu}}
\affil[1]{Department of Electrical Engineering and Computer Science, University of Tennessee, Knoxville, TN 37996, USA}
\maketitle



\begin{abstract}
 Koopman analysis can be used to understand the dynamics of a nonlinear dynamical system in terms a linear, but generally infinite dimensional operator.  The isostable coordinate system focuses on the slowest decaying principal Koopman eigenmodes.  This work leverages the isostable coordinate framework in the identification of slow manifolds for dynamical systems with fixed point attractors, defined as surfaces for which the fastest decaying isostable coordinates are zero.  Numerical challenges associated with separation between fast and slow timescales necessitate the development of new computational approaches to identify these slow manifolds.  Two such strategies are developed which approximate backward-time solutions on the slow manifold starting near the fixed point and extending far beyond the linear regime. Application to a variety of examples illustrates the utility of these methods and their potential use for model order reduction purposes.
\end{abstract}


\section{Introduction}

Model order reduction is often a critical first step for mathematical analysis and control design in engineering applications involving high-dimensional dynamical systems.  Projection-based methods \cite{benn15} such as proper orthogonal decomposition \cite{holm96}, \cite{town18}, \cite{berk93} and balanced truncation \cite{moor81}, \cite{guge04}, \cite{sore02} are particularly useful for systems that can be well-approximated through linearization.  By contrast, however, dynamical systems with nonnegligible nonlinearities usually require more targeted approaches that exploit specific features of the underlying vector fields.  

This work focuses on model order reduction of dynamical systems with a pronounced time scale separation between `fast' and `slow' dynamics.  Center manifold theory \cite{wigg03} represents one of the earliest examples of such a model order reduction framework, whereby the long term dynamics decay to a center manifold exponentially fast allowing the system dynamics to be well-approximated by the behavior on the center manifold. A variety of methods consider fast-slow subsystems that explicitly separate `fast variables' and `slow variables' of a dynamical system.  Singular perturbation theory can subsequently be used to understand the system dynamics in terms of the behavior on a slow manifold for which the system dynamics are strongly attracted to (resp.,~repelled from) provided the manifold is stable (resp.,~unstable) \cite{feni79}, \cite{kape99}.  These ideas have been extended in recent years to include slow manifolds of saddle-type stability \cite{farj18}, \cite{guck09}.  Recent work has also focused on understanding dynamical behavior on slow decaying invariant manifolds.  For instance, a spectral submanifold \cite{hall16} represents the smoothest invariant manifold that functions as a nonlinear extension of the associated linear modal subspace.  A reduced order model can be obtained by focusing on the spectral submanifold associated with slowest decaying solutions  \cite{pons20}, \cite{cene22}.  

In recent years, ideas stemming from Koopman operator theory have also been used for model order reduction purposes.  The Koopman operator can be used to represent a fully nonlinear dynamical system using a linear, but infinite dimensional operator \cite{budi12}, \cite{mezi13}, \cite{mezi19}.   Practical use of this theory often involves the search for Koopman invariant subspaces \cite{brun16}, \cite{lusc18}, \cite{kais21}, or in the case of dynamic mode decomposition, \cite{schm10}, \cite{will15}, \cite{kutz16} least-squares fitting of observables to linear models; both of these approaches typically require lifting the observable space to a {\it higher} dimension to capture relevant nonlinear behaviors.  By contrast, for dynamical systems with a stable attractor, the principal Koopman eigenmodes, i.e.,~the isostable coordinates, can be used to define an intrinsic reduced order coordinate system  \cite{maur13}, \cite{wils16isos}, \cite{wils17isored}. By focusing on the slowest decaying modes of the Koopman operator and ignoring the fast decaying terms, the isostable coordinate paradigm has been used for reduced order modeling in a variety of control applications \cite{maur16}, \cite{soot17}, \cite{wils21dd}, \cite{wils21input}, \cite{ahme23}.

This work focuses on the use of the isostable coordinate framework to aid in the definition and computation of a slow manifold for dynamical systems with a stable fixed point.  For a general dynamical system
\begin{equation} \label{maineq}
    \dot{x} = F(x),
\end{equation} 
where $x \in \mathbb{R}^N$ with a stable fixed point $x_0$ for which $F(x_0) = 0$, the isostable coordinate system can be used to characterize the decay of solutions evolving under the flow of the vector field \cite{maur13} in the basin of attraction of the fixed point.  Letting $\lambda_1,\dots,\lambda_N$ be the eigenvalues associated with the linearization of \eqref{maineq} and noting that the isostable coordinates represent level sets of Koopman eigenmodes, under the evolution of the flow of \eqref{maineq}  the principal isostable coordinates evolve according to  $\dot{\psi}_k = \lambda_k \psi_k$, where $\psi_k$ is an isostable coordinate with associated decay rate $\lambda_k$ \cite{wils20ddred}, \cite{wils21dd}.  Ordering the isostable coordinates in terms of their decay rate so that $| {\rm Real}(\lambda_k)| \leq |{\rm Real}(\lambda_{k+1})|$, provided there is some separation between $| {\rm Real}(\lambda_k)|$ and  $ |{\rm Real}(\lambda_{k+1})|$ a straightforward definition of a slow manifold is given by:
\begin{equation} \label{slowdef2}
    W^s = \{ x \in \mathbb{R}^N | \psi_k(x) = 0 \; {\rm for} \; k > \beta  \}.
\end{equation}
Because of the rapid decay of the faster isostable coordinates, for reduced order modeling purposes, it is typically possible to understand the dynamics of \eqref{maineq} in terms the behavior on the slow manifold.  Previous works \cite{wils15}, \cite{wils16isopde}, \cite{wils19cdc} have focused on obtaining reduced order representations for carefully chosen trajectories that approach the fixed point in infinite time, but not for the entire slow manifold as defined in \eqref{slowdef2}.  Other works consider asymptotic expansions in a basis of isostable coordinates \cite{maur16}, \cite{wils21dd}, but these approaches require the computation of high order partial derivatives and are generally computationally prohibitive to implement for high dimensional systems.

As a primary contribution, this work develops two general approaches for computing slow manifolds of dynamical systems defined according to Equation \eqref{slowdef2}.  These approaches involve the backward-time integration of solutions starting near the fixed point of \eqref{maineq} and extending in the basin of attraction far beyond the linear regime.  Due to the timescale separation inherent between the fast and slow isostable coordinates, direct backward-time integration is impossible, necessitating the development of the new techniques presented in this work.  The organization of this paper is as follows:~Section \ref{backsec} provides relevant background on isostable coordinates.  Section \ref{compsec} discusses two strategies for numerically computing the slow manifold; by directly computing the gradient of the slow isostable coordinates along trajectories and providing an approximation for the subspace spanned by the gradient of the fast decaying isostable coordinates, it is possible to overcome the aforementioned timescale separation problem.  Section \ref{ressec} shows results for this approach applied to a variety of dynamical systems and Section \ref{concsec} gives concluding remarks.

\section{Background on Isostable Coordinates} \label{backsec}

Consider the general differential equation given in \eqref{maineq}.  Letting $x_0$ be a stable fixed point with $F(x_0) = 0$, through local linearization one finds
\begin{equation} \label{lineq}
   \Delta \dot{x} = J \Delta x + O(||\Delta x||^2),
\end{equation}
where $\Delta x = x - x_0$ and $J$ is the Jacobian evaluated at $x_0$.  While Equation \eqref{lineq} itself is only useful in the limit that $\Delta x$ is small, its spectrum can be used to define a set of principal isostable coordinates, which represent level sets of principal Koopman eigenfunctions \cite{mezi19}, \cite{maur13}, \cite{mezi13} that correspond to states in the basin of attraction of the fixed point.  Specifically, let $w_k$, $v_k$ and $\lambda_k$ be left eigenvectors, right eigenvectors, and eigenvalues of $J$, respectively.  The eigenvalues will be ordered in terms of their decay rates so that $| {\rm Real}(\lambda_k)| \leq |{\rm Real}(\lambda_{k+1})|$.  For $\lambda_1$, the slowest decaying term, an associated principal isostable coordinated $\psi_1(x)$ can be defined by considering the infinite-time behavior of solutions
\begin{equation} \label{isodef}
    \psi_1(x) = \lim_{t \rightarrow \infty} (w_1^T (\phi(t,x) - x_0) \exp (- \lambda_1 t)),
\end{equation}
where $\phi(t,x)$ is the flow of \eqref{maineq} and $^T$ is the transpose.  Additionally, faster decaying principal isostable coordinates $\psi_2,\dots,\psi_N$ be defined implicitly as level sets of principal Koopman eigenfunctions with associated decay rates $\lambda_2, \dots, \lambda_N$ \cite{kval21}.    A hallmark of the isostable coordinates framework is that
\begin{equation} \label{psidyn}
    \dot{\psi}_k = \lambda_k \psi_k,
\end{equation}
for $k = 1,\dots,N$ when $x$ is in the basin of attraction of the fixed point evolving under the forward-time flow of \eqref{maineq}.  As a point of emphasis, while there are generally an infinite number of Koopman eigenfunctions (and hence an infinite number of isostable coordinates) the principal isostable coordinates $\psi_1,\dots,\psi_N$ are associated with the principal Koopman eigenfunctions (i.e.,~those associated with the linearization about $x_0$).  


In many cases (especially in situations where control is involved) it can be useful to compute the gradient of the isostable coordinate $I_k(x) =  \frac{\partial \psi_k}{\partial x}$ to characterize how trajectories are influenced by external perturbations.  Considering an arbitrary trajectory $x^\gamma(t) = \phi(t,x_0)$, $I_k$ can be computed straightforwardly along this trajectory according to the equation \cite{wils16isos}
\begin{equation} \label{isoeq}
\dot{I}_k = -(J^T - \lambda_k {\rm Id}) I_k,
\end{equation}
where $I_k$ is evaluated at $x^\gamma(t)$, $J$ is the Jacobian evaluated at $x^\gamma(t)$, and ${\rm Id}$ is an appropriately sized identity matrix.   A short proof of \eqref{isoeq} is included here for completeness:~considering a small, arbitrary perturbation $\Delta x$ to $x^\gamma(t)$ occurring at $t = 0$, let $x_\epsilon(t) = x^\gamma(t) + \epsilon \Delta x(t)$.  Under the flow of \eqref{maineq},
\begin{equation} \label{dxeq}
    \Delta \dot{x} = J \Delta x + O(||\Delta x||^2).  
\end{equation}
The change in $\psi_k$ resulting from this perturbation is 
\begin{equation} \label{dpsieq}
   \Delta \psi_k =  I_k^T \Delta x + O(||\Delta x||^2). 
\end{equation}
Taking the time derivative of \eqref{dpsieq} and truncating $O(||\Delta x||^2)$ terms yields  
\begin{align} \label{preadj}
    \Delta \dot{\psi}_k = \dot{I}_k^T \Delta x + I_k^T  \Delta \dot{x}.
\end{align}
Using the fact that $\Delta \dot{\psi}_k = \lambda_k \Delta \psi_k = \lambda_k I_k^T \Delta x$ and substituting \eqref{dxeq} into Equation \eqref{preadj} one finds
\begin{equation}
    0 = (\dot{I}_k^T + I_k^T(J - \lambda_k {\rm Id})) \Delta x.
\end{equation}
Noting that $\Delta x$ is arbitrary, Equation \eqref{isoeq} follows immediately.  In the specific case that $x^\gamma(t) = x_0$ for all $t > 0$, Equation \eqref{isoeq} simplifies to $0 = (J^T - \lambda {\rm Id}) I_k$ with solution $I_k = w_k$.  As such, any time that $x^\gamma(t)$ is close to $x_0$, \begin{equation} \label{ijeq}
I_k(t) = w_k + O(|| \Delta x||).
\end{equation}

\section{Definition and Computation of Slow Manifolds Using Isostable Coordinates} \label{compsec}

\subsection{Slow Manifolds of Linear Systems}  \label{linmodsec}
Consider linear model
\begin{equation} \label{linmod}
     \dot{x} = A  x,
\end{equation}
where $x \in \mathbb{R}^N$ and $A\in \mathbb{R}^{N \times N}$ has $N$ eigenvalues with negative real component sorted in order of their decay rates $| {\rm Real}(\lambda_k)| \leq |{\rm Real}(\lambda_{k+1})|$.  It will be assumed that $J$ is diagonalizable so that the the eigenvalues are not defective.  Let $v_j$ and $W_j$ be right and left eigenvalues associated with the eigenvalue $\lambda_j$.  Considering Equation \eqref{isodef}  the $\psi_1$ isostable coordinate of \eqref{linmod} is easily computed as $\psi_1(x) = w_1^T  x$.  Using methods from \cite{wils21dd}, it is straightforward to show that the faster decaying isostable coordinates for this linear system can be computed according to
\begin{equation} \label{linisos}
    \psi_k(x) = w_k^T x,
\end{equation}
with corresponding solution
\begin{equation}
    x(t) =  \sum_{j = 1}^N \psi_j(x(0)) \exp(\lambda_k t) v_j.
\end{equation}
A slow invariant manifold invariant manifold comprised of the slowest decaying contributions of the solution can be defined according to 
\begin{align} \label{linslow}
    W^s &= {\rm span}\{ v_1,\dots,v_\beta  \} \nonumber \\
    &= \{ x \in \mathbb{R}^N | \psi_k(x) = 0 \; {\rm for} \; k  > \beta  \},
\end{align}
with $\beta < N$.  When ${\rm Real}(\lambda_{\beta+1}) \ll {\rm Real}(\lambda_{\beta})$, solutions collapse rapidly to $W^s$.

\subsection{Slow Manifolds of Nonlinear Systems} 
Isostable coordinates can be used to extend the notion of a slow manifold from Section \ref{linmodsec} to a nonlinear system.  In a similar manner, considering a general system of the form \eqref{maineq}, one can define a slow manifold as
\begin{equation} \label{slowdef}
    W^s = \{ x \in \mathbb{R}^N | \psi_k(x) = 0 \; {\rm for} \; k > \beta  \}.
\end{equation}
Compared to the definition \eqref{linslow}, while the slow manifold is no longer a linear eigenspace, it can still be defined in terms of the principal isostable coordinates.  Considering the dynamics of isostable coordinates from \eqref{psidyn} it is immediately apparent that $W^s$ is an invariant manifold since $\dot{\psi}_k = 0$ for $k > \beta$ for all $x \in W^s$. 

\begin{figure}[htb]
\begin{center}
\includegraphics[height=2.0 in]{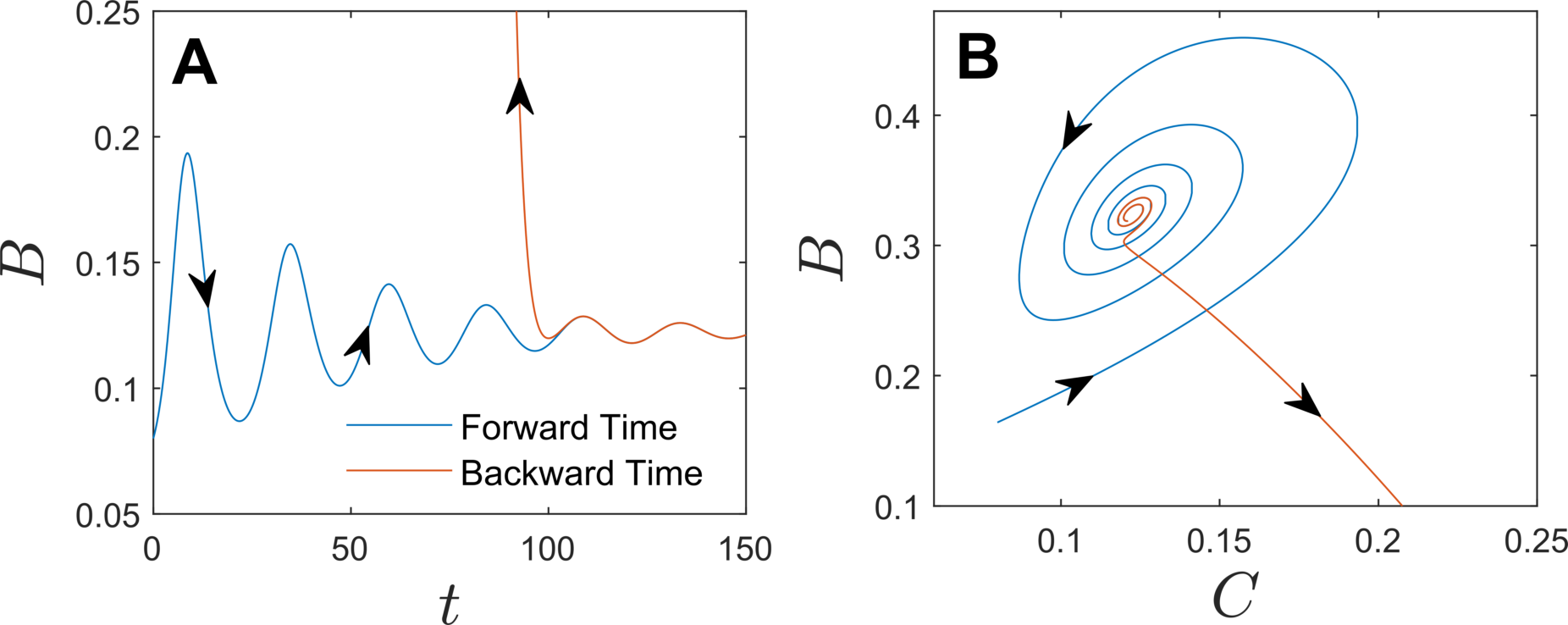}
\end{center}
\caption{Trajectories of \eqref{circmodel} display oscillations characteristic of a slow manifold associated with the complex conjugate eigenvalues $\lambda_{1,2}$ in forward time (blue lines).  This slow manifold can be estimated with accuracy to $O(||\Delta x||^2)$ near the fixed point, but simple backward integration cannot follow the slow manifold in backward time (red lines) to provide an estimate of the slow manifold extended to the nonlinear regime.  Panels A and B show trajectories from the same simulation.}
\label{backwardsets}
\end{figure}

The critical challenge that this work addresses is the computation of the slow manifold.    As a naive approach, If one could find some initial condition on $W^s$, forward and backward integration in time would trace out the profile of the manifold. 
 Considering \eqref{ijeq}, near the fixed point, $\psi_k(x) = w_k(x-x_0) + O(||\Delta x||^2)$.  Thus, $x = \sum_{j = 1}^\beta \alpha_j v_j$, where $\alpha_j = O(\epsilon)$ is $O(||\Delta x||^2)$ close to $W^s$.  However, backward time integration in time will quickly amplify the magnitude of $\psi_{\beta+1}, \dots, \psi_N$.  Even if it were possible to specify an initial condition exactly on the slow manifold, nearby initial conditions diverge in backward time at a rate governed by the fast time scales -- as such, inevitable errors from any ODE solver will be rapidly amplified.  As a concrete example, consider a three dimensional model based on the Goodwin oscillator taken from \cite{gonz05}
 \begin{align} \label{circmodel}
\dot{B} &= h_1 \frac{K_1^n}{K_1^n + D^n}  -  h_2  \frac{B}{K_2+B}  +  \alpha,   \nonumber \\
\dot{C} &= h_3 B - h_4 \frac{C}{K_4+C},  \nonumber \\
\dot{D} &=  h_5 C - h_6\frac{D}{K_6+D} , 
\end{align}
with parameters $n = 6$, $h_1 = 0.84$, $h_2 = 0.42$, $h_3 = 0.7$, $h_4 = 0.35$, $h_5 = 0.7$, $h_6 = 0.35$, $K_1 = K_2 = K_4 = K_6 = 1$, and $\alpha = 0.025$.  For this parameter set, the model has a stable fixed point at $(B,C,D) = (0.12,0.32,1.84)$ with eigenvalues $\lambda_{1,2} = -0.022 \pm 0.26$ and $\lambda_3 = -0.53$.  The separation between the eigenvalues make this model a good candidate for a reduction through the computation of a slow manifold defined according to Equation \eqref{slowdef} taking $\beta =2$, but it is not possible to follow the slow manifold in backward time for very long, as illustrated in Figure \ref{backwardsets}.  This numerical issue becomes more pronounced when the spectral gap between fastest and slowest decaying eigenvalues increases.

\subsection{Practical Numerical Computation of Trajectories Backwards in Time Along the Slow Manifold}

Considering any trajectory $x(t) \in \mathbb{R}^N$, the evolution of the isostable coordinate $\psi_k$ follows 
\begin{equation} \label{gradevol}
    \frac{d \psi_k}{d t} = \frac{\partial \psi_k}{\partial x}^T \left( \frac{dx}{dt} \right) = I_k^T \frac{dx}{dt}.
\end{equation}
Considering \eqref{psidyn} and letting $\tilde{t} = -t$ under the backward time flow of \eqref{maineq} $d \psi/ d\tilde{t} = -\lambda_k \psi_k$; combining this fact and \eqref{gradevol}, the evolution along the slow manifold as defined by \eqref{slowdef} in backward time is
\begin{equation}  \label{timeevolution}
    \frac{dx}{d \tilde{t}} = \begin{bmatrix} I_1^T \\ \vdots \\ I_\beta^T \\ I_{\beta+1}^T \\ \vdots \\ I_{N}^T   \end{bmatrix} ^{-1} \begin{bmatrix}  - \lambda_1 \psi_1 \\ \vdots \\ -\lambda_{\beta} \psi_{\beta} \\  0 \\ \vdots \\ 0 \end{bmatrix},
\end{equation}
where $^{-1}$ denotes the matrix inverse.  Provided  $I_1, \dots, I_N$ can be computed accurately and that the matrix inverse exists, the evolution in backward time along the slow manifold is given by \eqref{timeevolution}.  Different strategies will be necessary for computation of each $I_k$ depending on whether it corresponds to a fast or slow decaying isostable coordinate.

\subsection{Computation of $\boldsymbol{I_k}$ associated with Slow Decaying Isostable Coordinates} \label{slowdecay}

As shown below, $I_j$ can be accurately computed along the slow manifold in backwards time provided the spectral gap between $\lambda_1$ and $\lambda_j$ is small relative to the time scale of integration.   To begin, consider some $\hat{x}(t)$ evolving under the backward time flow of \eqref{maineq}, i.e.,~subject to $\dot{x} = -F(x)$ on the slow manifold defined in \eqref{slowdef} for some initial condition $\hat{x}(0)$ close to the fixed point.      Comparing with \eqref{isoeq}, under the evolution of the flow, i.e.,~letting $\tilde{t} = -t$ 
\begin{equation} \label{backtimeiso}
    \frac{ d {I}_j}{d \tilde{t}} = (J^T - \lambda_j {\rm Id})I_j.
\end{equation}
 Noticing that Equation \eqref{backtimeiso} is linear and time varying, let $\Phi_j(\tilde{t}_2,\tilde{t}_1)$ be the associated state transition matrix, i.e.,~with the property that $I_j(\tilde{t}_2) = \Phi_j(\tilde{t}_2,\tilde{t}_1) I_j(\tilde{t}_1)$.  With this in mind, noticing that for any $j$ and $k$
\begin{equation}
     \frac{ d {I}_j}{d \tilde{t}} = (J^T - \lambda_k {\rm Id}) I_j + (\lambda_k-\lambda_j) I_j,
\end{equation}
one can show
\begin{equation} \label{phisubs}
\Phi_j(\tilde{t}_2,\tilde{t}_1) = \Phi_k(\tilde{t}_2,\tilde{t}_1) \exp( (\lambda_k-\lambda_j) (\tilde{t}_2-\tilde{t}_1)  ).
\end{equation}

From \eqref{ijeq}, $w_j$ provides an $O(||\Delta x||)$ close estimate for $I_j$ at $\tilde{t} = 0$.  Letting $\hat{I}_k(\tilde{t} \, )$ denote the approximation of $I_k(\tilde{t} )$ obtained according to \eqref{backtimeiso}, at $\tilde{t} = 0$ using the initial condition $\hat{I}_j(0) = w_j$, 
\begin{equation}
    \hat{I}_j(0) = I_j(0) + \epsilon \sum_{k = 1}^N s_k(0) I_k(0),
\end{equation}
where the truncated $O(||\Delta x||)$ terms are written in a basis comprised of $I_1(0), \dots, I_N(0)$ with $s_k(0)$ being the associated coordinates in this basis.   Letting $h$ be a constant time step of a numerical solver, the solution after the first time step at $\tilde{t} = h$ is 
\begin{align}
    \hat{I}_j(h) &= \Phi_j(h,0)\left( I_j(0) + \epsilon \sum_{k = 1}^N s_k(0) I_k(0) \right) +  \epsilon \sum_{k = 1}^N s_k(h) I_k(h) \nonumber \\ 
    &= I_j(h) + \epsilon \sum_{k=1}^N s_k(0) I_k(h)(\exp((\lambda_k - \lambda_j)h) +  \epsilon \sum_{k = 1}^N s_k(h) I_k(h).
    \end{align}
Above, an additional $O(\epsilon)$ error resulting from the numerical solver is represented in a basis comprised of $I_1(0), \dots, I_N(0)$ and the second line is obtained after substituting \eqref{phisubs} and simplifying.  Repeating this process over $M$ time steps, one finds
\begin{align} \label{errorseq}
    \hat{I}_j(Mh) &= I_j(Mh) + \epsilon \sum_{k = 1}^N s_k(0) I_k(Mh) \exp((\lambda_k-\lambda_j)Mh) \nonumber \\
    & + \epsilon \sum_{i = 1}^M  \sum_{k = 1}^N s_k(ih) I_k(Mh) \exp((\lambda_k - \lambda_j)((M-i)h)).
\end{align}
For small enough $h$, the final term in \eqref{errorseq} can be approximated by an integral leaving
\begin{align} \label{errorapprox}
    \hat{I}_j(\tilde{t}) &\approx   I_j(\tilde{t}) + \epsilon \sum_{k = 1}^N s_k(0) I_k(\tilde{t}) \exp((\lambda_k-\lambda_j)\tilde{t}) \nonumber \\
    & \quad + \frac{\epsilon}{h}  \sum_{k = 1}^N \int_0^{\tilde{t}} s_k(\tau) I_k(\tilde{t}) \exp((\lambda_k-\lambda_j)(\tilde{t} -\tau)) d\tau,
\end{align}
where $\tilde{t} = Mh$ is the length of integration.  Considering Equation \eqref{errorapprox},  focusing on the error $E_j(\tilde{t}) = \hat{I}_j(\tilde{t})-I_j(\tilde{t})$ and letting $\bar{s}_k = \max_{h} | s_k(i h)|$
\begin{align} \label{errorbound}
||E_j(\tilde{t})|| &    \leq    \epsilon \sum_{k = 1}^N \bar{s}_k || I_k(\tilde{t})|| \exp((\lambda_k - \lambda_j) \tilde{t}) +   \frac{\epsilon   }{h}  \sum_{k = 1}^N \bar{s}_k ||I_k(\tilde{t})||   \int_0^{\tilde{t}}   \exp((\lambda_k-\lambda_j)(\tilde{t} -\tau)) d\tau \nonumber \\
&=   \epsilon \sum_{k = 1}^N \bar{s}_k || I_k(\tilde{t})|| \exp((\lambda_k - \lambda_j) \tilde{t}) +  \frac{\epsilon   }{h}  \sum_{k = 1}^N  \bar{s}_k ||I_k(\tilde{t})|| \bigg( \frac{\exp((\lambda_k-\lambda_j)\tilde{t}) - 1}{\lambda_k - \lambda_j} \bigg).
\end{align}
Considering the terms of \eqref{errorbound}, numerical errors made at each time step in practical numerical solvers are typically of the form $K h^\alpha$ where $\alpha\geq 2$ \cite{asch98} and $K$ is an $O(1)$ term.  Given that the error at each time step is assumed to be $O(\epsilon)$, $\frac{1}{h} = O(1/\sqrt[\alpha]{\epsilon})$.  It will be assumed that $\bar{s}_k ||I_k||$ is smaller than $O(1/\epsilon)$.  As such, the errors in the computation of $I(\tilde{t})$ remain small provided that $\exp((\lambda_k - \lambda_j)\tilde{t})$ remains an $O(1)$ term for all $k$.  In general, this allows for the accurate computation of $I_1,\dots,I_\beta$ (i.e.,~the terms associated with the slowly decaying isostable coordinates) when the spectral gap between the eigenvalues $|{\rm Real}(\lambda_1 - \lambda_\beta)|$ is small.  In general, however, the errors in the computation $I_{\beta+1},\dots,I_N$ will grow rapidly necessitating computation using a different strategy.


\subsection{Alternative to Direct Computation of $\boldsymbol{I_k}$ Associated with Fast Decaying Isostable Coordinates}

Section \ref{slowdecay} illustrates that \eqref{backtimeiso} can only be used to numerically compute $I_1, \dots, I_\beta$ associated with the slow isostable coordinates on the slow manifold accurately.  Considering \eqref{timeevolution}, along the slow manifold, $dx/d\tilde{t}$ is orthogonal to  ${\rm Span}(W_{\beta+1}^N)$ where $W_{\beta+1}^N = \{ I_{\beta+1},\dots,I_N \}$.  With this in mind, it is only necessary to find a basis for ${\rm Span}(W_{\beta+1}^N)$ in order to implement \eqref{timeevolution}.

Towards this goal, let $g_k \in \mathbb{C}^N$ be defined so that 
\begin{equation} \label{isonorm}
    I_j^T g_k = \begin{cases} 1, & \text{if } k = j, \\
    0, & \text{otherwise}.
    \end{cases}
\end{equation}
From the definition above, letting $G_1^\beta = \{ g_1,\dots,g_\beta \}$, ${\rm Span}(G_1^\beta)^\bot = {\rm Span}(W_{\beta+1}^N)$, shifting the focus towards finding a basis for ${\rm Span}(G_1^\beta)$ whose orthogonal complement provides the remaining information about ${\rm Span}(W_{\beta+1}^N)$.

To proceed, one again considering an arbitrary trajectory $x^\gamma(t) = \phi(t,x_0)$, $g_k$ evolves along this trajectory according to 
\begin{equation} \label{gdyn}
    \dot{g}_k = (J - \lambda_k {\rm Id})g_k,
\end{equation}
where $g_k$ is evaluated at $x^\gamma(t)$ and $J$ is the Jacobian evaluated at $x^\gamma(t)$.  This can be verified using both \eqref{isoeq} and \eqref{gdyn} to show that $\frac{d}{dt} \big( I_k^T g_j \big) = 0$ for any $k$ and $j$.  Once again, under the backward time flow letting $\tilde{t} = -t$, Equation \eqref{gdyn} becomes
\begin{equation} \label{gback}
    \frac{d g_k}{d \tilde{t}} = -(J + \lambda_k {\rm Id}) g_k.
\end{equation}
Noticing that \eqref{gback} is the adjoint system of \eqref{backtimeiso}, it has the solution $g_k(\tilde{t}_2) = \Phi_k^T(\tilde{t}_1,\tilde{t}_2) g_k(\tilde{t}_1)$, where $\Phi_k$ is the state transition matrix associated with the solution of \eqref{backtimeiso} and $^T$ denotes the transpose.  Recalling in Section \ref{slowdecay} that $\Phi_k (\tilde{t}_2,\tilde{t}_1)$ attenuated numerical errors for solutions evaluated along backwards trajectories, $\Phi_k (\tilde{t}_1,\tilde{t}_2)$ (and hence $\Phi_k^T(\tilde{t}_1,\tilde{t}_2)$)  will amplify numerical errors, ultimately precluding numerical computation of $g_1,\dots,g_\beta$ using \eqref{gback}. Alternative strategies to compute ${\rm Span}(G_1^\beta)$ are discussed below.


\subsection{Computation of Trajectories Along the Slow Manifold Using Asymptotic Expansions for    \boldmath{$g_1,\dots,g_\beta$}   }     \label{asymethod}

The terms $g_k$ for $k = 1,\dots,\beta$ can be approximated by first noting from \eqref{psidyn}, that for solutions evolving under the forward time flow of \eqref{maineq},
\begin{equation} \label{psidot2}
\dot{\psi}_k = \lambda_k \psi_k = I_k^T \frac{dx}{dt}.
\end{equation}
Considering \eqref{isonorm} and \eqref{psidot2} together, this implies
\begin{equation} \label{gderiv}
    \frac{dx}{dt} = \sum_{j = 1}^N \lambda_j \psi_j g_j
\end{equation}
must be satisfied for any $x$ evolving under the flow of the vector field \eqref{maineq}.  Previous work \cite{wils21dd} considered the state $x$ in terms of a Taylor expansion in powers of the isostable coordinates
\begin{equation} \label{xexp}
    x = x_0 + \sum_{k = 1}^N \psi_k v_k + \underbrace{\sum_{j = 1}^N \sum_{k = 1}^j \psi_j \psi_k h^{jk}}_{ \text{second order terms}} +   \underbrace{\sum_{i = 1}^N \sum_{j = 1}^i \sum_{k = 1}^j \psi_i \psi_j \psi_k h^{ijk}}_{ \text{third order terms}} + \dots ,
\end{equation}
where $v_k$ is an eigenvector associated with the eigenvalue $\lambda_k$ and $h^{jk}, h^{ijk}, \dots$ $\in \mathbb{C}^N$.  Using a strategy described in \cite{wils21dd}, it is straightforward to compute each  $h^{jk}, h^{ijk}, \dots$ with a description given in Appendix \ref{apxb}.  Taking the time derivative of \eqref{xexp} yields
\begin{equation} \label{deriveq}
      \frac{dx}{dt} =  \sum_{k = 1}^N \lambda_k \psi_k v_k + \sum_{j = 1}^N \sum_{k = 1}^j (\lambda_k + \lambda_j)\psi_j \psi_k h^{jk} +   \sum_{i = 1}^N  \sum_{j = 1}^i \sum_{k = 1}^j (\lambda_k + \lambda_j) \psi_i \psi_j \psi_k h^{jk} + \dots.
\end{equation}
Subsequently, equating \eqref{deriveq} and \eqref{gderiv} and collecting terms associated with $g_1$, one finds
\begin{equation} \label{ord3}
    g_1 = v_1 + \underbrace{\bigg( 2 h^{11}\psi_1 + \sum_{k \neq 1} \psi_k h^{k1} \bigg)}_{ \text{second order terms}} + \underbrace{\bigg( 3  \psi_1^2 h^{111} + 2 \sum_{k \neq 1}  \psi_1 \psi_k h^{k11} + \sum_{k \neq 1} \sum_{i \neq 1} \psi_k \psi_i h^{ij1} \bigg)}_{\text{third order terms}} + \dots.
\end{equation}
Similar relationships can be used to compute $g_2,\dots,g_\beta$ to arbitrary accuracy in the expansion in the basis of isostable coordinates.  This information can be used to approximate ${\rm Span}(W_{\beta+1}^N)$ by finding the span of $(G_1^\beta)^\bot$.  subsequently, trajectories can be computed backwards in time starting from some initial condition close to the fixed point while  computing $\frac{dx}{dt}$ according to Equation \eqref{timeevolution} and concurrently computing $I_1,\dots,I_\beta$ according to \eqref{backtimeiso}.  This accuracy of this method method of computation is limited by the accuracy of the asymptotic expansion of $g_1,\dots,g_\beta$.

\subsection{An Alternative Predictor-Corrector Strategy for Computation of Trajectories along the Slow Manifold} \label{pcmethod}

The strategy discussed in Section \ref{asymethod} requires the computation of the terms $v_k, h^{jk}, h^{ijk}, \dots$ from the asymptotic expansion \eqref{deriveq},  As discussed in Appendix \ref{apxa}, this can be prohibitively difficult for high dimensional systems.  An alternative predictor-corrector approach for computation of trajectories along the slow manifold is discussed here.  This strategy still uses \eqref{timeevolution} to compute backwards time trajectories and uses \eqref{backtimeiso} to compute $I_1,\dots,I_\beta$.  The key difference in this strategy is in the approximation of $g_1,\dots,g_\beta$ and the subsequent determination of ${\rm Span}(W_{\beta+1}^N)$.

\subsubsection{Approximation of \boldmath{$g_1,\dots,g_\beta$}}
For this approach, it will be assumed that that $|\lambda_k| + |\lambda_j| <  |\lambda_\beta|$ for all $k,j \leq \beta$ so that Equation \eqref{phijeq} derived in Appendix \ref{apxa} can be used to approximate the state transition matrix, $\Phi_J(t_2,t_1)$ for the linear time-varying equation 
\begin{equation} \label{jaceq}
    \Delta \dot{x} = J(t) \Delta x,
\end{equation} 
where $J$ is the Jacobian of \eqref{maineq} evaluated for any trajectory $x(t)$ that approaches $x_0$ in the limit that $t$ goes to infinity.

To proceed, consider Equation \eqref{gdyn}, which gives the dynamics of $g_k$ for solutions $x(t)$ evolving under the flow of \eqref{maineq}.  Noticing that Equation \eqref{gdyn} is a linear time varying equation, its solution is given by $g_k(t_2) = \Phi_{g_k}(t_2,t_1) g_k(t_1)$ where $\Phi_{g_k}$ is the state transition matrix.  Additionally, noting the similarity between Equations \eqref{jaceq} and \eqref{gdyn}, one can show that their state transition matrices are related by
\begin{align}
\Phi_{g_j}(t_2,t_1) &= \Phi_J(t_2,t_1) \exp(-\lambda_j (t_2-t_1)).
\end{align}
Substituting \eqref{phijeq} into the above equation, provided $t_2-t_1$ is large enough, $\Phi_{g_j}(t_2,t_1)$ is well approximated by
\begin{align} \label{phigapprox}
\Phi_{g_j}(t_2,t_1) \approx \sum_{k = 1}^\beta  \hat{v}_{k}  \hat{\lambda}_{k}\exp(-\lambda_j (t_2-t_1)) \hat{w}_{k}^T  +   \sum_{k = \beta  +1}^N  \hat{v}_{k}  \hat{\lambda}_{k}\exp(-\lambda_j (t_2-t_1)) \hat{w}_{k}^T,
\end{align}
where $\mathcal{V}_\beta = {\rm Span}(\{\hat{v}_1,\dots,\hat{v}_\beta\}) = {\rm Span}(\{ v_1,\dots,v_\beta \})$ and $\hat{\lambda}_k = O(\epsilon)$ for $k>\beta$.  When $t$ is large, $x(t)$ is close enough to the fixed point so that $g_j(t_2) = v_j + O(\epsilon)$.  Therefore,
\begin{align} \label{gjapprox}
    g_j(t_1) &= \Phi^{-1}_{g_j}(t_2,t_1) (v_k + O(\epsilon)) \nonumber \\
    & \approx   \sum_{k = 1}^\beta   \frac{1}{\hat{\lambda}_{k}\exp(-\lambda_j (t_2-t_1))}  \hat{v}_{k}  \hat{w}_{k}^T (v_j + O(\epsilon))  +   \sum_{k = \beta  +1}^N  \frac{1}{  \hat{\lambda}_{k}\exp(-\lambda_j (t_2-t_1)) } \hat{v}_{k}  \hat{w}_{k}^T (v_j + O(\epsilon)).
\end{align}
Above, the second line follows after substituting Equation \eqref{phigapprox}.  
\begin{figure}[htb]
\begin{center}
\includegraphics[height=1.8 in]{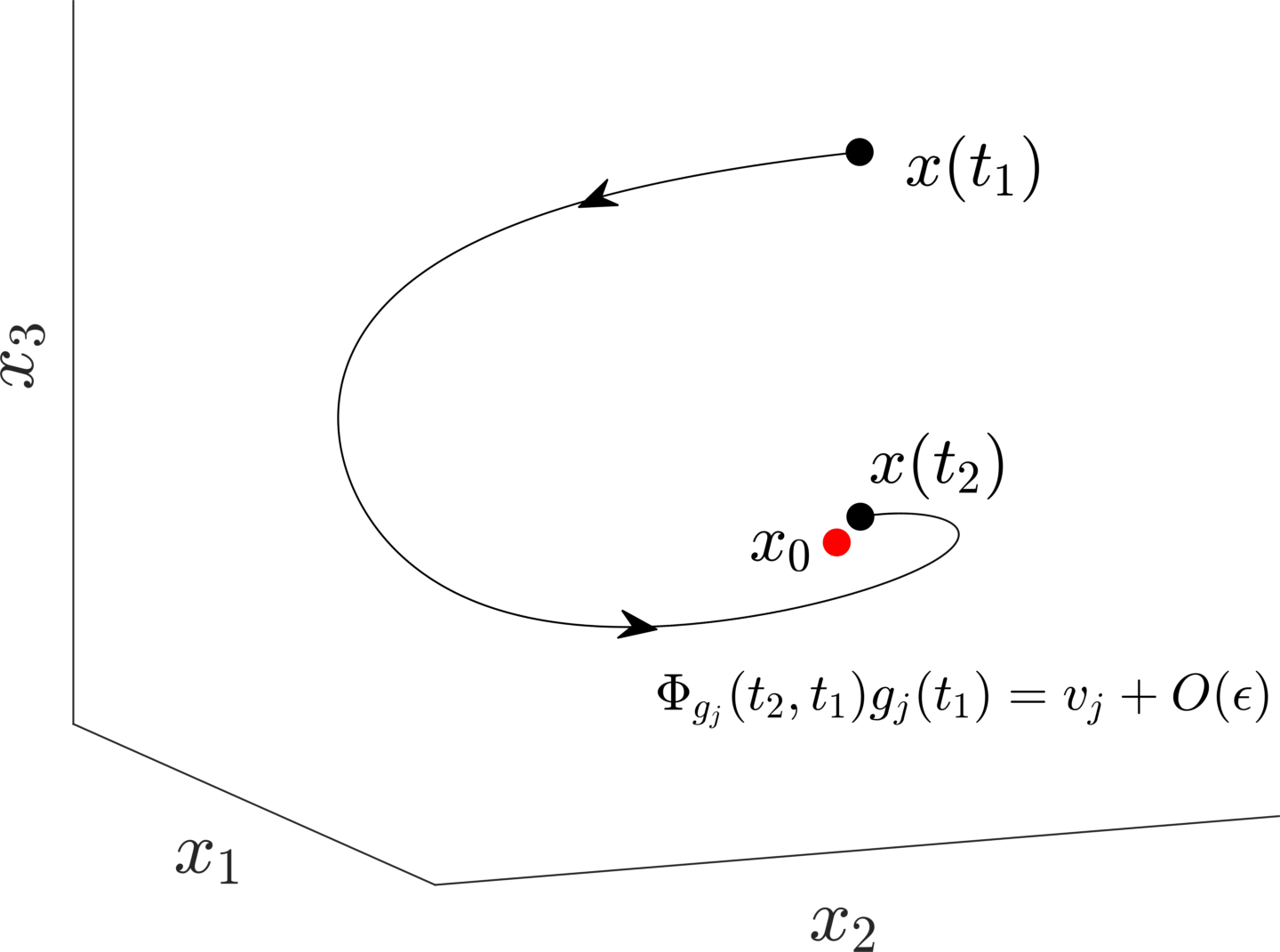}
\end{center}
\caption{The predictor-corrector strategy detailed below requires an accurate approximation of $g_j(t_1)$.  While $g_j(t_1)$ can be computed according to $g_j(t_1) = \Phi^{-1}_{g_j}(t_2,t_1) (v_k + O(\epsilon))$, this is not practically feasible when some of the eigenvalues of $\Phi_{g_j}(t_2,t_1)$ are near zero.  Equation \eqref{gjapprox2} provides a computationally tractable approximation which is valid when $g_k(t_1)$ is $O(\epsilon)$ close to a linear subspace spanned by $\{{v}_1,\dots,{v}_\beta\}$. }
\label{examppc}
\end{figure}

In Equation \eqref{gjapprox}, notice that $1/ ( \hat{\lambda}_{k}\exp(-\lambda_j (t_2-t_1)) )$ is an $O(1/\epsilon)$ term for $k = \beta+1,\dots,N$ which makes direct computation numerically infeasible.  As an alternative approach, consider the matrix $P_\beta = \begin{bmatrix} \hat{v}_1&\dots& \hat{v}_{\beta}\end{bmatrix} \in \mathbb{C}^{N \times \beta}$ and the relationship
\begin{equation}
   g_{j}^\perp(t_1) = P_\beta P_\beta^\dagger g_j(t_1),
\end{equation}
where $\dagger$ is the pseudoinverse provides the projection of  $g_j(t_1)$ onto $\mathcal{V}_\beta$.  The residual, i.e.,~$g_j(t_1)-g_j^\perp(t_1)$ is then
\begin{align} \label{rjeq}
    r_j(t_1) &= ({\rm Id} - P_\beta P_\beta^\dagger) g_j(t_1) \nonumber \\
     &= \begin{bmatrix} \hat{v}_{\beta+1} & \dots & \hat{v}_{N}  \end{bmatrix} \alpha,
\end{align}
where $\alpha \in \mathbb{C}^{N-\beta}$ with the $k^{\rm th}$ entry being equal to $(1/  \hat{\lambda}_{k+\beta}\exp(-\lambda_j (t_2-t_1))) \hat{v}_{k+\beta} \hat{w}_{k+\beta}^T (v_j + O(\epsilon)) $ per Equation \eqref{gjapprox}.  It will be assumed that $||r_j(t_1)||_2 = O(\epsilon)$ and $\sigma_{\rm min} \big(  \begin{bmatrix} \hat{v}_{\beta+1} & \dots & \hat{v}_{N}  \end{bmatrix} \big) = O(1)$ where $|| \cdot ||_2$ denotes the two norm and $\sigma_{\rm min}$ gives the minimum singular value.  Provided that these two assumptions are valid it follows immediately from Equation \eqref{rjeq} that $||\alpha||_2 = O(\epsilon)$.  Subsequently, Equation \eqref{gjapprox} can be rewritten truncating $O(\epsilon)$ terms as
\begin{equation} \label{gjapprox2}
    g_j(t_1) \approx     \sum_{k = 1}^\beta   \frac{1}{\hat{\lambda}_{k}\exp(-\lambda_j (t_2-t_1))}  \hat{v}_{k}  \hat{w}_{k}^T v_j.
\end{equation}
Of course, the accuracy of \eqref{gjapprox2} is contingent on the validity of the assumptions  $||r_j(t)||_2 = O(\epsilon)$ and $\sigma_{\rm min} \big( \begin{bmatrix} \hat{v}_{\beta+1} & \dots & \hat{v}_{N}  \end{bmatrix} \big) = O(1)$.  The first assumption holds when $x(t_1)$  is near the fixed point $x_0$ so that $g_j(t_1) = v_j + O(\epsilon)$, and will generally begin to degrade as the state moves farther from the fixed point. The second assumption can be checked directly with knowledge of $\Phi_{g_j}(t_2,t_1)$.  The condition \eqref{gjapprox2} can be used as part of a prediction-correction method as described below to compute trajectories backwards in time along the slow manifold.

\subsubsection{Prediction Step}  \label{predsec}

Choose some initial condition with $x(t_2) = \sum_{k = 1}^\beta  \alpha_k v_k$ where $\alpha_1,\dots,\alpha_\beta$ are small so that $\psi_i(t_2) = \alpha_i + O(\epsilon)$.  Integrate along the slow manifold backwards starting at time $t_2$ and ending at time $t_1$.  For this computation $\frac{dx}{dt}$ can be computed according to Equation \eqref{timeevolution} using ${\rm Span}(W_{\beta+1}^N) = (G_1^\beta)^\bot$ where $G_1^\beta$ is approximated by  $\{ v_1 \; \dots \; v_\beta   \}$; this approximation will introduce some error in the computation that will be mitigated during the correction step.  The terms $I_1,\dots,I_\beta$ can be computed on the interval $t \in [t_2,t_1]$ according to \eqref{backtimeiso}.  Also compute $\Phi_{g_j}(t_2,t_1)$. 

\subsubsection{Correction Step} \label{correctsec}
For the prediction step, the approximation ${\rm Span}(W_{\beta+1}^N) = (G_1^\beta)^\bot$ with $G_1^\beta \approx \{ v_1 \; \dots \; v_\beta   \}$ introduces some error into the computation of $x(t)$.  As such, $x_2$ will not be exactly on the slow manifold.  This error can be mitigated by identifying the correction  $\Delta x$ that places the $x(t_1)$ on the slow manifold.  To proceed, considering \eqref{gderiv}
\begin{align} \label{resid_0}
\frac{dx}{dt} = \sum_{j = 1}^\beta \lambda_j \psi_j(t_1) g_j(t_1),
\end{align}
 where each $g_j(t_1)$ can be approximated according to Equation \eqref{gjapprox2}.  Additionally, from \eqref{maineq} 
\begin{align} \label{resid_00}
\frac{dx}{dt} &= F(x(t_1) + \Delta x) \nonumber \\
&= F(x(t_1)) + J \Delta x,
\end{align}
where $J$ is the Jacobian of \eqref{maineq} evaluated at $x(t_1)$ and $\Delta x$ is the correction to $x(t_1)$.  Equating \eqref{resid_0} and \eqref{resid_00} yields
\begin{equation} \label{resid1}
    \sum_{j = 1}^\beta \lambda_j \psi_j(t_1) g_j(t_1) = F(x(t_1)) + J \Delta x.
\end{equation}
Perturbations in the span of $\{ \hat{v}_{\beta+1}, \dots,  \hat{v}_{\beta+N}  \}$ will not influence the isostable coordinates $\psi_1, \dots \psi_\beta$.  With this in mind, the $\Delta x$ which minimizes the residual between the left and right sides of \eqref{resid1} and does not change the isostable coordinate is given by
\begin{equation} \label{eqncorr}
    \Delta x = \begin{bmatrix} \hat{v}_{\beta+1} & \dots & \hat{v}_{N}  \end{bmatrix} \begin{bmatrix} \hat{v}_{\beta+1} & \dots & \hat{v}_{N}  \end{bmatrix}^\dagger  J^\dagger \bigg(\sum_{j = 1}^\beta \lambda_j \psi_j(t_1) g_j(t_1) - F(x(t_1)) \bigg).
\end{equation}
Using the corrected forward time initial condition $x(t_1) + \Delta x$,  $\Phi(t_2,t_1)$ and $x(t)$ can be recomputed in forward time by directly simulating Equation \eqref{maineq}.  Using the resulting $x(t)$, a correction for $I_1(t_1),\dots,I_\beta(t_1)$ can be obtained through computation of \eqref{isoeq} in backwards time.  

\subsubsection{Implementation of the Predictor-Correction Strategy Over Multiple Iterations}
Suppose that the predictor-corrector strategy described in Sections \eqref{predsec} and \eqref{correctsec} has been implemented at least once, yielding an an approximation for $I_1(t),\dots,I_\beta(t)$, $x(t)$ and $\Phi(t_2,t)$ for $t \in [t_1,t_2]$ where $t_1 < t_2$.  Using the prediction step, taking initial conditions $x(t_1)$ and $I_1(t_1),\dots,I_\beta(t_1)$, an estimate for $I_1(t),\dots,I_\beta(t)$ and $x(t)$ can be obtained for $t \in [t_1 -\Delta t,t_1]$ where $\Delta t > 0$.  This information can also be used to compute $\Phi(t_1,t_1 - \Delta t)$ as well as $\Phi(t_2,t_1 - \Delta t) = \Phi(t_2,t_1)\Phi(t_1,t_1 - \Delta t)$.  Next, a correction for $x(t_1 - \Delta t)$ can be computed using \eqref{eqncorr} along with a subsequent approximation of $x(t)$, $I_1(t), \dots I_\beta(t)$, and  $\Phi(t_2,t)$ on the interval $t \in [t_2,t_1 - \Delta t]$.  This procedure can be repeated over multiple iterations.  Smaller values of $\Delta t$ will limit the error that is allowed to accumulate before the correction step, however, additional correction steps will require additional computational time.  $\Delta t$ must be chosen to balance this trade-off. 

As mentioned earlier, this predictor-corrector strategy will generally break down as $||r_j(t)||_2$ from \eqref{rjeq} becomes larger.  While this condition cannot generally be checked explicitly, the examples to follow will illustrate that the predictor-corrector strategy provides a more accurate approximation for the slow manifold than the strategy from \eqref{asymethod} based on asymptotic expansion.

\section{Results} \label{ressec}

The numerical strategies from Sections \ref{asymethod} and \ref{pcmethod} for computing trajectories along a slow manifold defined by \eqref{slowdef} are illustrated below for a variety of dynamical systems.

\subsection{Slow Manifold of a Planar Model} \label{slowplan}

\begin{figure}[htb]
\begin{center}
\includegraphics[height=2.5 in]{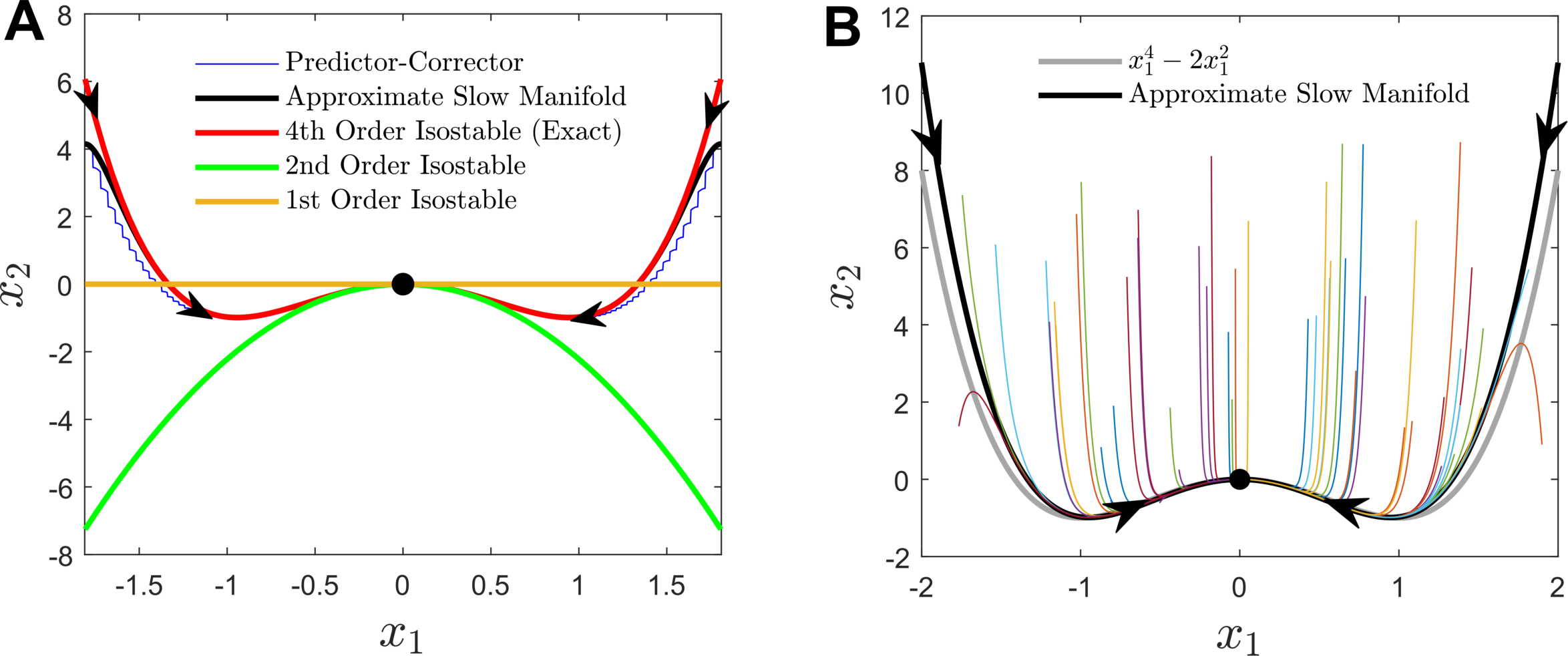}
\end{center}
\caption{Orange, green, and red curves in Panel A approximations of the slow manifold using the strategy from \ref{asymethod} where $g_1$ is computed to different orders of accuracy.  Higher order approximation of $g_1(\psi_1)$ yields a more accurate estimate of the slow manifold.  For this example, the fourth order estimate of $g_1(\psi_1)$ is exact.  The blue line shows the result of the predictor-corrector method strategy described in Section \ref{pcmethod}.  Integrating the result forward in time gives a close approximation of the true slow manifold (black line).  Panel B plots an approximation of the slow manifold using the predictor-corrector method (black line) against the $\dot{x}_2$ nullcline (grey line).  Solutions of \eqref{planmod}, shown as colored lines, rapidly converge to the slow manifold.}
\label{twodimmodelresult}
\end{figure}

Consider the following planar system:
\begin{align} \label{planmod}
    \dot{x}_1 &= -0.05 x_1, \nonumber \\
    \dot{x}_2 &= -(x_2 - x_1^4 + 2 x_1^2).
\end{align}
Solutions of \eqref{planmod} slowly decay towards the fixed point at $x_1=x_2=0$ closely following the $\dot{x}_2 = 0$ nullcline $x_2 = x_1^4 - 2x_1^2$.  This general behavior is shown in Panel B of Figure \ref{twodimmodelresult}.  Local linearization yields eigenvalues $\lambda_1 = -0.05$ and $\lambda_2 = -1$ with associated eigenvectors $v_1 = [1\;0]^T$ and $v_2  = [0\;1]^T$.   For this example, a one-dimensional slow manifold will be computed for which $\psi_2 = 0$.  Starting with an initial condition $x =  0.001 v_1$ at $t_2 = 0$ with corresponding isostable coordinates $\psi_1 \approx   0.001$ and $\psi_2 \approx 0$, both methods described in Sections \ref{asymethod} and \ref{pcmethod} are used to follow the slow manifold backwards in time.  This process is repeated for an initial condition $x =  -0.001 v_1$ to compute the portion of the unstable manifold with negative isostable coordinates. 

Results are shown in Figure \ref{twodimmodelresult}.  When using the strategy from Section \ref{asymethod}, the asymptotic approximation of $g_1(\psi_1)$ from Equation \eqref{ord3} is taken to first, second, and fourth order accuracy in the isostable coordinates; the resulting approximations of the slow manifold for which $\psi_2 = 0$ is shown in orange, green, and red, respectively, in panel A.  The first order approximation  for $g_1(\psi)$ is a straight line.  Increasing the order of accuracy of the approximation of $g_1(\psi_1)$ yields a better approximation of the actual slow manifold

Note here that the fifth order terms and beyond are all zero so that the fourth order expansion provides an exact representation of $g_1(\psi)$ subsequently allowing for the slow manifold to be computed with negligible errors.  The blue line in panel A shows the results of the predictor-corrector approach using $\Delta t = 0.5$ for each prediction step with a total of 300 steps.  Using the predictor-corrector approach, an approximation of the slow manifold (black line) is obtained through forward integration of the state obtained after the final prediction step at $t_1 = 150$.  While the black and red lines do not match perfectly the differences quickly become negligible as the system converges towards the true slow manifold.  

In panel B, the predictor-corrector method is extended for an additional 10 steps and the resulting approximation of the slow manifold (black line) is plotted against the $\dot{x}_2 = 0$ isocline (grey line).  Colored lines show the evolution of 50 randomly chosen initial conditions evolving under the flow of \eqref{planmod}, with states rapidly converging to the approximation for the slow manifold.   For the simple planar model equations \eqref{planmod}, the strategy from Section \ref{asymethod} is the best choice for computation of the slow manifold because the expansion for $g_1(\psi_1)$ from \eqref{ord3} is exact when taken to fourth order accuracy.  However, generally the asymptotic expansion will not be exact at any order of accuracy.

\subsection{Pendulum with Variable Damping} \label{pendsec}
Consider the equations of motion for a pendulum with variable damping:
\begin{align} \label{pendvar}
    \dot{x}_1 &= x_2, \nonumber \\
    \dot{x}_2 &= -\alpha x_1 - \sin(x_1) - (\beta + \gamma x_3)x_2, \nonumber \\
    \dot{x}_3 &= - \kappa (x_3 - r^2),
\end{align}
where $x_1$ and $x_2$ represent the angular position and angular velocity, respectively, of the pendulum, $r^2 = x_1^2 + x_2^2$, $\alpha = 0.23$ is a torsional spring constant, and $(\beta + \gamma x_3)$ sets a viscous damping coefficient where $\beta = 0.1$ and $\gamma = 0.01$.  $x_3$ is an auxiliary variable with that contributes to the variable damping.  Taking $\kappa = 8$, the value of $x_3$ quickly converges to the surface $x_3 = x_1^2  +x_2^2$. The model \eqref{pendvar} has a globally stable fixed point at $x_1 = x_2 = x_3 =0$.  Local linearization yields eigenvalues $\lambda_{1,2} = -0.05 \pm 1.11 i$ and $\lambda_3 = -8$.   For this example, a two-dimensional slow manifold will be computed along which $\psi_3 = 0$.

\begin{figure}[htb]
\begin{center}
\includegraphics[height=2.5 in]{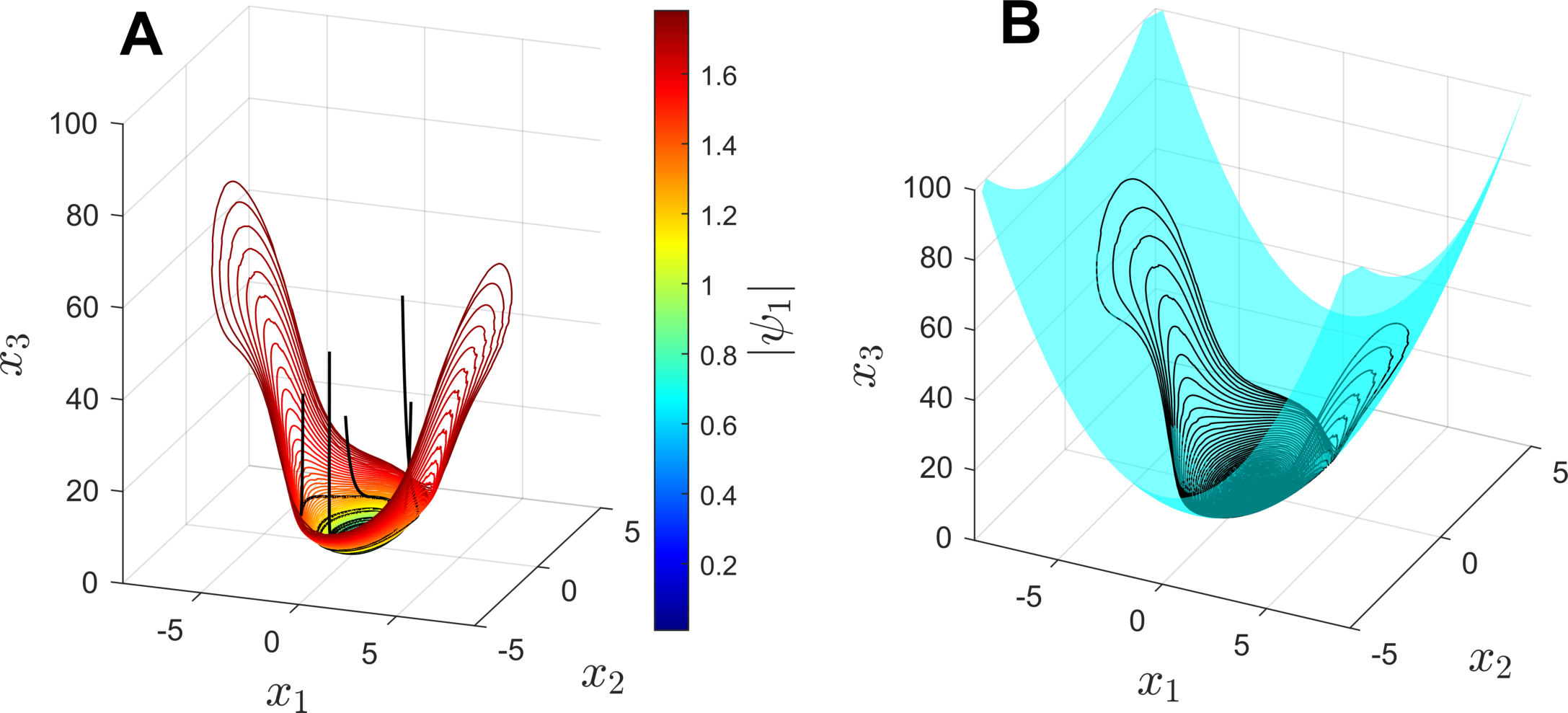}
\end{center}
\caption{Colored lines in panel A show level sets of $|\psi_1|$ on the slow manifold of \eqref{pendvar} computed according to the predictor-corrector strategy.  Example trajectories in black rapidly decay to to the slow manifold before converging to the fixed point. Panel B shows the same level sets plotted against the $\dot{x}_3 = 0$ nullcline,  $x_3 = x_1^2$.}
\label{pendlevelsets}
\end{figure}

Figure \ref{pendlevelsets} shows the slow manifold computed according using the predictor-corrector strategy from Section \ref{pcmethod}.  To do so, a set of initial conditions with $\psi_1 \approx 0.01 \exp(i k/200)$, $\psi_2 \approx 0.01 \exp(-i k/200)$ for $k = 1,\dots,200$ and $\psi_3 \approx 0$ are integrated backwards in time along the slow manifold.  Colored lines in panel A show level sets of $|\psi_1|$ computed according to this strategy on the slow manifold.  Black lines show five different initial conditions evolving under the flow of \eqref{pendvar} that rapidly converge to the slow manifold on the way to the fixed point.  Panel B shows the same curves for the level sets of $|\psi_1|$ on the slow manifold plotted against the surface $x_3 = x_1^2  +x_2^2$.  As expected, the $\dot{x}_3 = 0$ level set lies very close to the slow manifold.  

\begin{figure}[htb]
\begin{center}
\includegraphics[height=2.5 in]{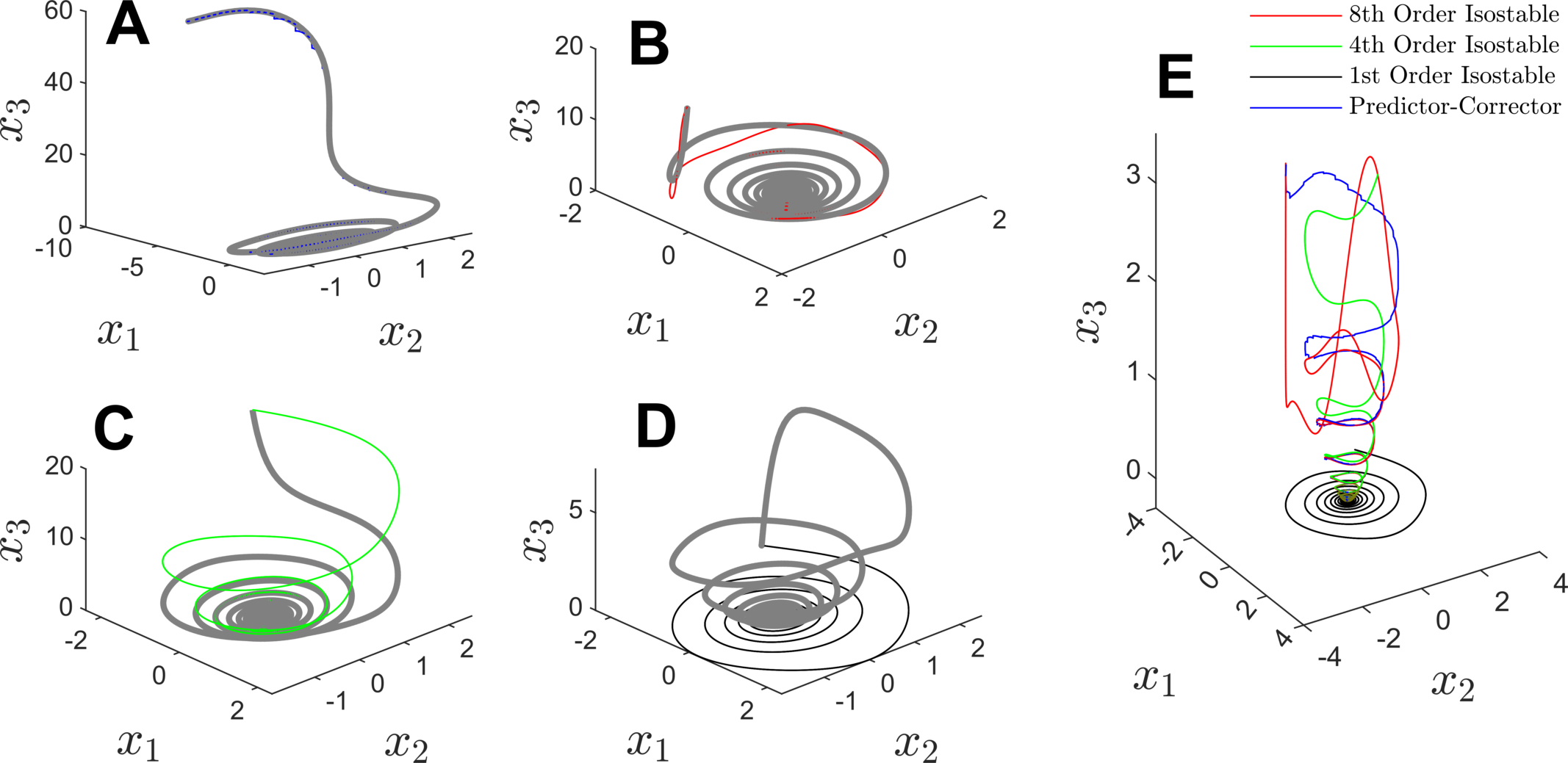}
\end{center}
\caption{Computation of trajectories along the slow manifold using the predictor-corrector method from Section \ref{pcmethod} (panel A) and the asymptotic approximation strategy \ref{asymethod} (panels B, C, and D). Colored lines in each panel show the approximation of the solution along the unstable manifold and the grey line shows the resulting trajectory starting at an initial condition corresponding to the largest isostable coordinate and integrating \eqref{pendvar} forwards in time.   Panel E shows each trajectory from panels $A-D$ plotted on the same axis.    }
\label{individualsolutions}
\end{figure}

Figure \ref{individualsolutions} shows the result of using both methods described in Sections \ref{asymethod} and \ref{pcmethod} to compute the backwards time trajectory along the slow manifold with an initial condition $\psi_{1,2} \approx 0.01$ and $\psi_3 \approx 0$.  The approach from Section \ref{asymethod} approximates $g_1$ and $g_2$ with an asymptotic expansion in a basis of isostable coordinates and is not sufficient to compute this trajectory as the state moves farther from the fixed point. Panel A of Figure \ref{individualsolutions} shows the approximation for the predictor-corrector method computed 100 time units in backward time (blue line). Note that the discontinuities come from the correction step which occurs every 0.25 time units.  The grey line shows the resulting trajectory computed in forward time; the agreement between the forward and backward time curves indicates that this is a good approximation for the slow manifold.  In panels B, C, and D, the red, green, and black lines are computed taking the asymptotic expansion of $g_1$ up to eighth, fourth, and first order accuracy, respectively. Grey lines in each of these plots show the resulting trajectory starting at an initial condition corresponding to the largest isostable coordinate and integrating \eqref{pendvar} forward in time.  Noting that the slow manifold is invariant under the flow of the vector field, mismatch between the grey and colored curves in panels B, C, and D indicates that none of these solutions starts on the slow manifold.  Panel E shows each resulting approximations plotted on the same axes.  As the order of accuracy in the asymptotic expansion of $g_1$ increases, the resulting solution follows the true solution (well-approximate by the blue curve) for longer.  Computing the expansion of $g_1$ past eighth order accuracy becomes computationally infeasible.

\FloatBarrier

\subsection{Reduced Order Modeling of Nonlinear Oscillations in the Goodwin Model} \label{circsec}
\begin{figure}[htb]
\begin{center}
\includegraphics[height=2.5 in]{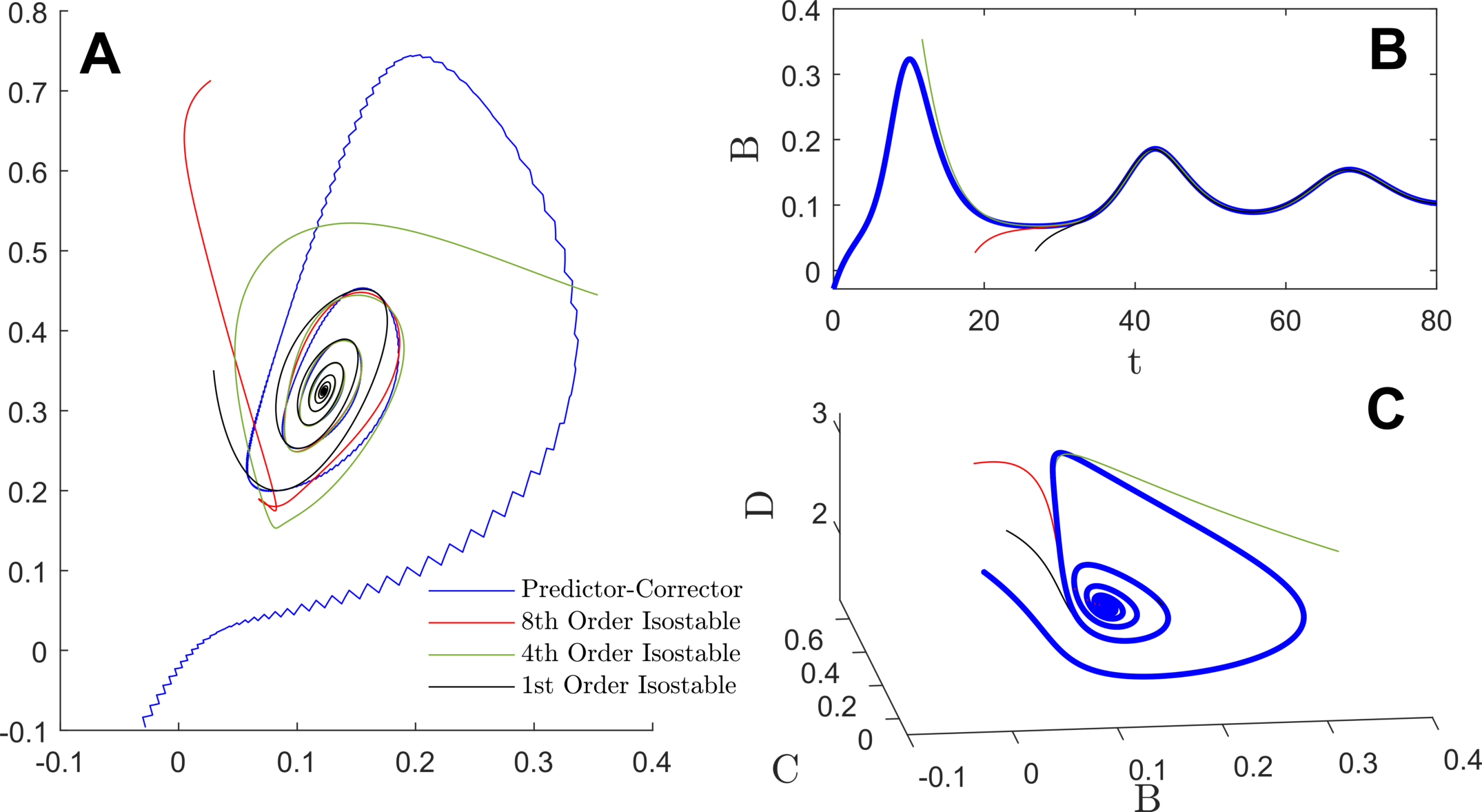}
\end{center}
\caption{For the Goodwin oscillator from Equation \eqref{circmodel} panel A shows trajectories computed along the slow manifold using the predictor-corrector method from Section \ref{pcmethod} and the asymptotic approximation strategy \ref{asymethod}. Subsequently integrating forwards in time, traces of the resulting trajectories are plotted in panels B and C.  In general, the predictor-corrector approach can be used to accurately compute the trajectory along the slow manifold further from the fixed point than the strategies based on asymptotic expansion of $g_1$.}
\label{circsolutions}
\end{figure}

The slow manifold is computed for the three-dimensional model introduced in Equation \eqref{circmodel}. The fixed point at $(B,C,D) = (0.12,0.32,1.84)$ has eigenvalues $\lambda_{1,2} = -0.022 \pm 0.26$ and $\lambda_3 = -0.53$.  For this example, a two-dimensional slow manifold will be computed along which $\psi_3 = 0$.  Figure \ref{circsolutions} shows the result of using both methods described in Sections \ref{asymethod} and \ref{pcmethod} to compute the backwards time trajectory along the slow manifold with an initial condition $\psi_{1,2} \approx 0.001$ and $\psi_3 \approx 0$.  In panel A, the blue line shows the approximation of the slow manifold using the predictor-corrector approach from Section \ref{pcmethod}. The solution is integrated backward in time for 275 time units with a correction step is applied every 0.25 time units.  The solution is subsequently integrated forward in time and the resulting trajectory is plotted in panels B and C as a blue line.  The strategy from Section \ref{asymethod} is also used to compute the backwards time trajectory along the slow manifold taking the asymptotic expansion of $g_1$ up to eighth, fourth, and first order accuracy.  Note that when using the strategy from Section \ref{asymethod}, the backward integration does not return once it begins to diverge from the slow manifold.  For each of these curves, the solution is subsequently simulated forwards in time and plotted in panels B and C.  The prediction-corrector methods is able to accurately follow the solution backwards in time along the slow manifold much farther than the strategy based on asymptotic expansion of $g_1$. As with the previous example, computing the expansion of $g_1$ past eighth order accuracy becomes computationally infeasible.

\begin{figure}[htb]
\begin{center}
\includegraphics[height=2.8 in]{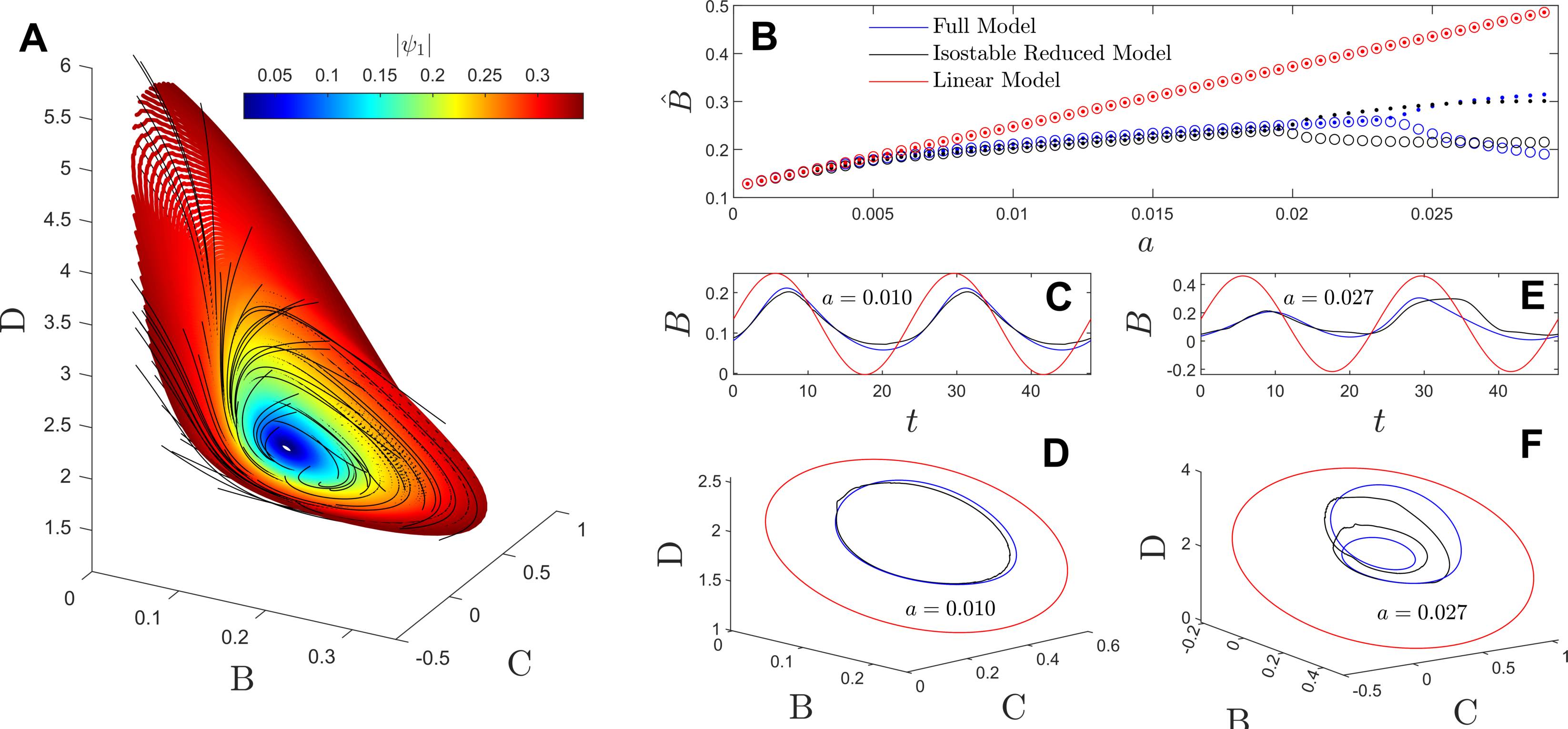}
\end{center}
\caption{Colored lines in panel A show level sets of $|\psi_1|$ on the slow manifold of \eqref{circmodel} computed according to the predictor-corrector strategy.  Example trajectories in black rapidly decay to to the slow manifold before converging to the fixed point.  The slow manifold can be used to capture the system response to an external input.  In panels B-F, an isostable-coordinate-based reduced order model of the form \eqref{redmodel} is compared to simulations of \eqref{circmodelinput} using the external input $u(t) = a \sin(2\pi t/T)$.  Panel $B$ shows the local maximum value of $B$ achieved over two successive cycles in steady state.  The period doubling bifurcation is captured with the isostable reduced model but not with the linearized model. Panels C-F show the steady state response to the sinusoidal input with different values of $a$.  Here, the isostable-based model is significantly more accurate than the linearized model.}
\label{circforcedresults}
\end{figure}

Figure \ref{circforcedresults} shows the slow manifold computed according using the predictor-corrector strategy from Section \ref{pcmethod}.  To do so, a set of initial conditions with $\psi_1 \approx 0.01 \exp(i k/200), \psi_2 \approx 0.01 \exp(-i k/200)$ for $k = 1,\dots,200$ and $\psi_3 \approx 0$ are integrated backwards in time along the slow manifold.  Colored lines in panel A of Figure \ref{circforcedresults} show level sets of $|\psi_1|$ computed according to this strategy which trace out the slow manifold.  Black lines show 100 different initial conditions evolving under the flow of \eqref{pendvar} that rapidly converge to the slow manifold on the way to the fixed point. 

After computing an approximation for the 2-dimensional slow manifold, it is straightforward to use this information for reduced order modeling purposes.  As an illustration, consider the Goodwin oscillator model from \eqref{circmodel} with the addition of an external input
 \begin{align} \label{circmodelinput}
\dot{B} &= h_1 \frac{K_1^n}{K_1^n + D^n}  -  h_2  \frac{B}{K_2+B}  +  \alpha + u(t),   \nonumber \\
\dot{C} &= h_3 B - h_4 \frac{C}{K_4+C},  \nonumber \\
\dot{D} &=  h_5 C - h_6\frac{D}{K_6+D} , 
\end{align}
where $u(t)$ is a time dependent input and all other terms and parameters are identical to those from \eqref{circmodel}.  When computing the trajectories that comprise the slow manifold according to the predictor-corrector method, the gradient of the isostable coordinate $I_1(\psi_1)$ is already calculated using Equation \eqref{isoeq}.  Changing variables to isostable coordinates, yields
\begin{align} \label{redmodel}
    \dot{\psi}_1 &= \frac{\partial \psi_1}{\partial x} \cdot \frac{dx}{dt} \nonumber \\
    &= \lambda_1 \psi_1  + i(\psi_1) u(t),
\end{align}
where $i(\psi_1) \equiv \begin{bmatrix} 1  & 0 & 0  \end{bmatrix} I_1(\psi_1)$.  Above, the second line can be obtained by noting that $\dot{\psi} = \lambda_1 \psi_1$ when $u(t) = 0$ and that $\frac{\partial \psi_1}{\partial x} = I_1(\psi_1)$ is solely a function of $\psi_1$ on the slow manifold  (with $\psi_2$ being the complex conjugate of $\psi_1$).  The output from the reduced order model $\begin{bmatrix} B(\psi_1) & C(\psi_1) & D(\psi_1) \end{bmatrix}$ is also solely a function of $\psi_1$ on the slow manifold.  Both $I(\psi_1)$ and $x(\psi_1)$ can be computed through interpolation of the individual trajectories used to compute the slow manifold.

The reduced order model \eqref{redmodel} is simulated with inputs $u(t) = a \sin(2 \pi t/24)$ for various values of $a$.  This model is simulated until a steady state response results for each value of $a$.  The simulations are compared to the full order model \eqref{circmodelinput} and to a local linearization about the stable fixed point with results shown in panels B-F of Figure \ref{circforcedresults}.  Open circles and dots in panel B show the local maximum value of $B$ achieved reached over two successive cycles in steady state.  In the full order and isostable reduced order models, a period doubling bifurcation occurs at $a = 0.023$ and $a = 0.020$, respectively. The linearized model does not display the period doubling bifurcation.  Panels C and D (resp.,~E and F) display the steady state solution taking $a = 0.010$ (resp.,~$a = 0.027)$.  Note that qualitatively similar results are also observed for sinusoidal inputs with other frequencies.

\FloatBarrier

\subsection{Reduced order Modeling of a Coupled Population of Planar Oscillators} \label{planarpop}

As a final example, consider a population of coupled planar oscillators
\begin{align} \label{isoclockadd}
\dot{x}_j &= \sigma x_j (\mu - r_j^2) - y_j (1+\rho_j (r_j^2 - \mu)) + \frac{K}{N}\sum_{i \neq j} x_i + u(t), \nonumber \\
\dot{y}_j &= \sigma y_j (\mu - r_j^2) + x_j (1+\rho_j ( r_j^2 - \mu)),
\end{align}
for $j = 1,\dots,N$ where $x_j$ and $y_j$ are the state variables of the $j^{\rm th}$ oscillator, $r_j^2 = (x_j^2 + y_j^2)$, and $N = 10$.  Here, $K = 1.54$ is the coupling strength, $u(t)$ is a common input applied to each oscillator, with constants $\mu = -4.5$, $\sigma = 0.05$ and $\rho_j = -0.2 + 4j/90$.  When $u(t) = 0$, Equation \eqref{isoclockadd} has a stable fixed point at $x_j = y_j = 0$ for all $j$.   Linearizing about this fixed point the slowest decaying eigenvalues are $\lambda_{1,2} =   -0.012 \pm 0.369i$.  The next slowest decaying eigenvalues are  $\lambda_{3,4} = -0.216 + 0.383i$.  In this example a two dimensional slow manifold will be computed for which $\psi_3 = \dots = \psi_{20} = 0$.

\begin{figure}[htb]
\begin{center}
\includegraphics[height=3.0 in]{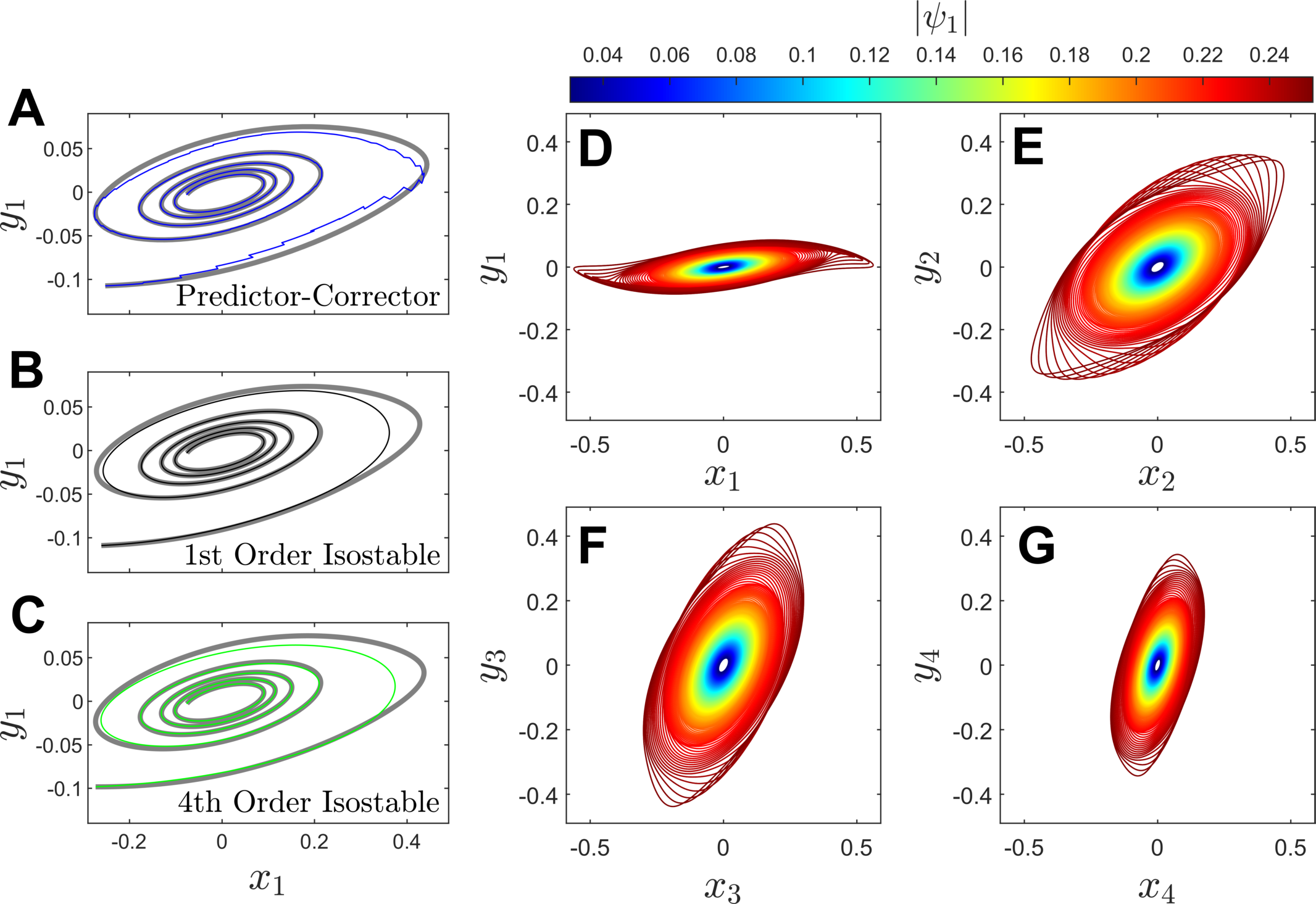}
\end{center}
\caption{Thin lines in Panel A (resp.,~B and C) show trajectories computed backwards in time along the slow manifold using the predictor-corrector (resp.,~asymptotic approximation) method.  The grey line shows the final condition obtained from these respective approximations plotted forward in time.  Panels D-G show level sets of $|\psi_1|$ on the slow manifold of \eqref{isoclockadd} for the first four oscillators computed according to the predictor-corrector strategy.}
\label{hopffigure}
\end{figure}

Panels A-C of  \ref{hopffigure} show the result of using both methods from Sections \ref{asymethod} and \ref{pcmethod} to compute the backward time trajectory along the slow manifold (taking $u(t) = 0$) with an initial condition $\psi_{1,2} \approx 0.1$ and $\psi_3 \approx \dots = \psi_{20} \approx 0$.  The blue line in Panel A is the approximation of the slow manifold obtained from the predictor-corrector approach from Section \ref{pcmethod}.  The solution is integrated backward in time for 84 time units with a correction step applied every 0.5 time units.  The result is then integrated forward in time and the resulting trajectory is plotted as a grey line.  Panels B and C show results using the strategy from Section \ref{asymethod} to compute the backwards time trajectory along the slow manifold taking the asymptotic expansion of $g_1$ to first and fourth order accuracy, respectively.  For this example, each of these methods yield a similar approximation for the slow manifold.  Note that the trajectory obtained using the first order approximation for $g_1$ is not simply a linear approximation of the slow manifold.   Panels D-G of Figure \ref{hopffigure} show level sets of $|\psi_1|$ on the slow manifold of \eqref{isoclockadd} computed according using the predictor-corrector strategy from Section \ref{pcmethod}.  To obtain these level sets, a set of initial conditions with $\psi_1 \approx 0.01 \exp(i k/200), \psi_2 \approx 0.01 \exp(-i k/200)$ for $k = 1,\dots,200$ and $\psi_3 \approx \psi_{20}  \approx 0$ is integrated backwards in time along the slow manifold.

\begin{figure}[htb]
\begin{center}
\includegraphics[height=2.5 in]{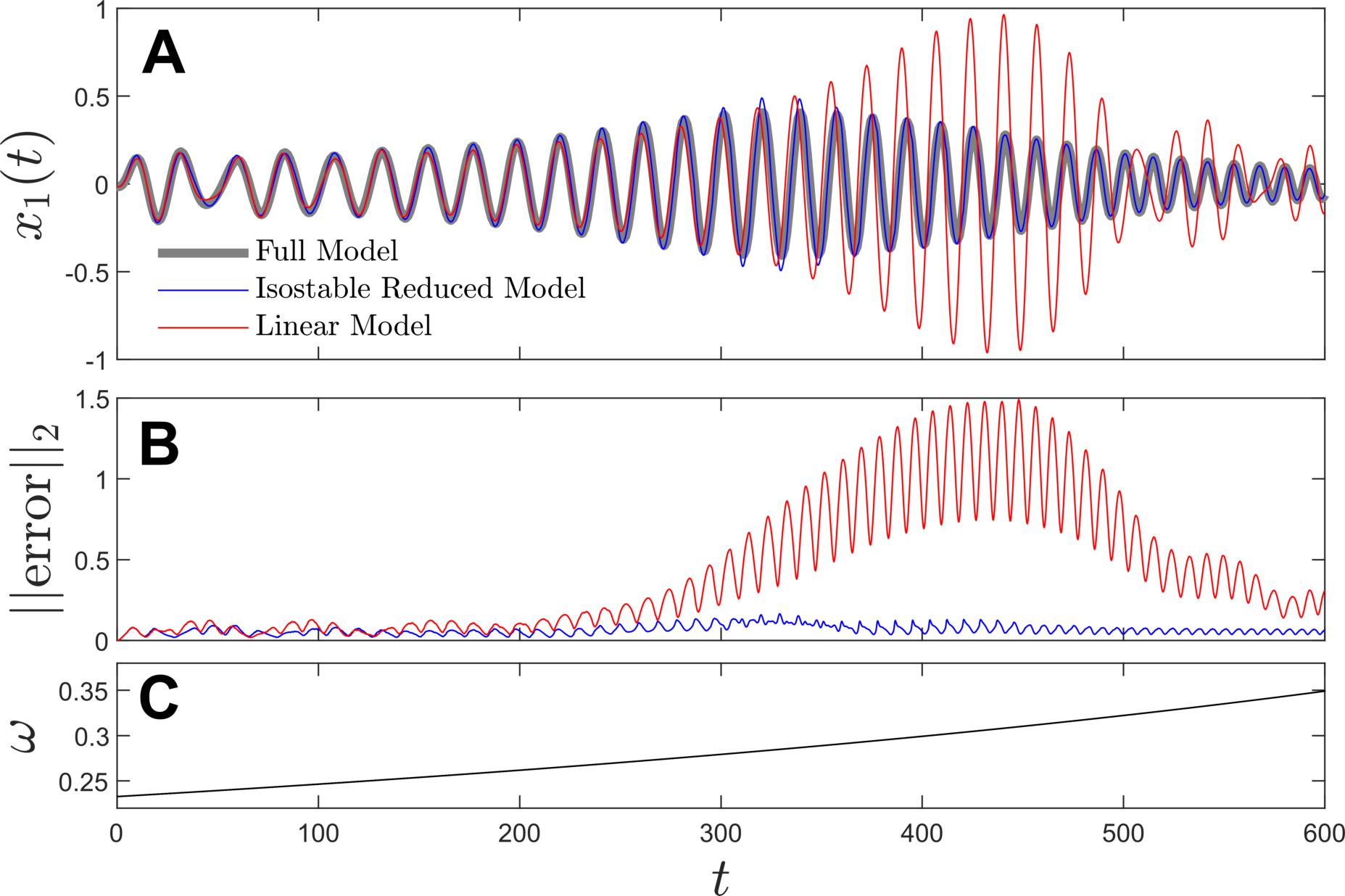}
\end{center}
\caption{Simulation of \eqref{isoclockadd}, the isostable coordinate reduced model \eqref{redhopf}, and a local linearization about the fixed point in response to the input $u(t) = 0.0045 \sin(\omega t)$. Panel A shows a trace of $x_1(t)$ for each model and panel B shows the 2-norm of the error between the isostable reduced model and linear model compared to the full model.  The variable frequency, $\omega$, of the input is plotted in panel C.  The forced behavior of the full order model is well-approximated by the dynamics on the slow manifold.}
\label{reducedsimresults}
\end{figure}

As with the example from Section \ref{circsec}, an isostable-coordinate-reduced order model is considered to capture the dynamics on the slow manifold of \eqref{isoclockadd}.  Mirroring the strategy that yields Equation \eqref{redmodel}, changing to isostable coordinates yields
\begin{equation} \label{redhopf}
    \dot{\psi}_1 = \lambda_1 \psi_1 + i(\psi_1) u(t),
\end{equation}
where $i(\psi_1) \equiv \begin{bmatrix} 1 & 0 & 1 & 0 & \dots & 1 & 0 \end{bmatrix} I_1(\psi_1)$, where $I(\psi_1)$ can be computed by interpolation of the individual trajectories computed on the slow manifold.  As a representative example, the reduced order model \eqref{redhopf} is simulated with input $u(t) = 0.0045\sin(\omega t)$; results are shown in Figure \ref{reducedsimresults} taking $\omega(t) = \frac{2 \pi}{27-0.15t}$, corresponding to a gradual increase in the forcing frequency.  This frequency band is chosen in order to capture the model response to forcing near its resonant frequency.  Panel A of Figure \ref{reducedsimresults} shows a trace of $x_1(t)$ over the course of this simulation, comparing the full model simulations \eqref{isoclockadd}, reduced order model simulations \eqref{redhopf}, and a model obtained from local linearization about the stable fixed point.  Panel B shows the 2-norm of the error (considering the states of all oscillators) for the same simulation.  The applied frequency is shown as a function of time in panel C.  The resonant peak of both the isostable reduced and full order model is observed when $\omega \approx 0.28$; the resonant peak of the linear model occurs closer to $\omega = 0.31$ and has a much larger amplitude than the resonant peak of the full order model.  Overall, the error associated with the isostable reduced model is substantially smaller than the error for the linear model, especially during the latter half of the simulation where the applied frequency is near the resonant peak. Additionally, the linear model displays spurious resonances when $\omega \in [0.32, 0.35]$ (in the final 100 time units of the simulation) that are not observed in the isostable reduced model.

\FloatBarrier

\section{Discussion and Conclusion}  \label{concsec}

This work considers the identification and computation of slow manifolds of dynamical systems with stable fixed points.   Provided the spectrum of the Koopman operator is discrete, the principal Koopman eigenfunctions can be used to define principal isostable coordinates, $\psi_1,\dots,\psi_N$.  Provided there is a large enough spectral gap between the decay rates of $\psi_\beta$ and $\psi_{\beta+1}$, a slow manifold can be defined as a surface for which $\psi_k = 0$ for $k > \beta$.  Numerical computation becomes challenging due to the spectral gap between fast and slow decaying isostable coordinates necessitating the development of two new strategies detailed in Section \ref{compsec}.  The examples considered in this work illustrate that the resulting dynamics in response to perturbations can be well-approximated by considering the behavior on the slow manifold.

While the general computational strategies detailed in this work represent a promising reduced order modeling approach, there are still many questions left to address.  Foremost, both of the computation strategies presented in Section \ref{compsec} require technical conditions to be satisfied that are difficult to explicitly check.  Specifically, the strategy suggested in Section \ref{asymethod} is valid provided the approximations of $g_1,\dots,g_\beta$ remain accurate to a given order of accuracy.  The predictor-corrector strategy from Section \ref{pcmethod} assumes that each $g_1,\dots,g_\beta$ can be well approximated by a projection onto a linear subspace spanned by $v_1, \dots, v_\beta$ (i.e.,~the slowly decaying eigenvectors associated with the linearization of \eqref{maineq}).  While both of of these assumptions are satisfied for states in a close vicinity of the fixed point in the linear regime,  there is no obvious metric to gauge how far these conditions extend into the nonlinear regime.  Despite this issue, the proposed computational strategies are able to provide approximations of the slow manifold that are still accurate for model order reduction purposes -- even if these computational strategies introduce moderately sized errors, forward time trajectories will converge exponentially fast to the slow manifold owing to the exponentially fast decay of the fast isostable coordinates.  

In principal, the strategies presented here can be applied to dynamical systems of arbitrarily high dimension.  For practical application, the predictor-corrector method from section \ref{pcmethod} is better suited to systems of exceptionally high dimension since it only requires numerical approximations of the Jacobian evaluated locally along a given trajectory.  By contrast, the strategy from Section \ref{asymethod} requires the symbolic computation of partial derivatives to orders that increase with the desired order of accuracy of $g_1,\dots,g_\beta$.  For the 20-dimensional system considered in Section \ref{planarpop}, this approximation can only be computed to fourth order accuracy before the computational burden becomes too large.  

 In this work, approximations of slow manifolds are obtained by choosing initial conditions near the fixed points and integrating backwards in time. For 1-dimensional manifolds, it is straightforward to implement this strategy with two initial conditions ($\psi_1 = \pm \epsilon$ where $0 < \epsilon \ll 1$).  For a two dimensional manifold with complex-conjugate isostable coordinates, it is similarly straightforward to choose a set of initial conditions for which $\psi_1 = \epsilon \exp(i \theta)$ for $\theta \in [0, 2\pi)$.  It is possible to apply this strategy to compute higher dimensional slow manifolds; in general an $n$-dimensional slow manifold will require the integration of an $n-1$ dimensional set of initial conditions backward in time along the slow manifold. Care would need to be taken in the choice of these initial conditions to obtain an accurate sampling on the slow manifold.

 It would be of general interest to consider extensions of the strategies developed in this work for dynamical systems with limit cycle attractors.  Previous work has considered the asymptotic expansion of the state for oscillatory systems in a basis of isostable coordinates \cite{wils20highacc}, but these strategies become computationally infeasible for high dimensional systems.  Additionally, it would be of interest to develop data-driven model identification strategies to infer the dynamics on the slow manifold in applications where the underlying model equations are unknown.  Additional questions of uniqueness of the resulting slow manifolds may arise in certain situations.  The expansion \eqref{xexp} is only unique provided appropriate nonresonance conditions on the eigenvalues of the linearization are satisfied \cite{wils20highacc}.  These considerations will ultimately influence on the uniqueness of the resulting slow manifold.  Similar nonresonance conditions on eigenvalues are required when considering the uniqueness of spectral submanifolds \cite{hall16}, \cite{pons20}, \cite{cene22}.  More careful investigation of these questions will be considered in future work.


This material is based upon the work supported by the National Science Foundation (NSF) under Grant No.~CMMI-2140527.

\begin{appendices}

\section{Perturbed Solutions of Trajectories Approaching the Fixed Point} \label{apxa}
\renewcommand{\thetable}{A\arabic{table}}  
\renewcommand{\thefigure}{A\arabic{figure}} 
\renewcommand{\theequation}{A\arabic{equation}} 
\setcounter{equation}{0}
\setcounter{figure}{0}



Consider Equation \eqref{maineq} with an initial condition $x(t_1)$ and corresponding isostable coordinates $\psi_1(t_1), \dots, \psi_N(t_1)$.  Consider also a small perturbation $x_1 = x(t_1) + \Delta x$ with corresponding isostable coordinates $\psi_1(t_1) + \Delta \psi_1, \dots, \psi_N(t_1) + \Delta  \psi_N$.   Using the Taylor expansion from Equation \eqref{xexp} difference $ \Delta x(t_1) = x(t_1)-x_1(t_1) = \Delta x$ can be written as
\begin{equation}
    \Delta x(t_1) = \sum_{k = 1}^N \Delta \psi_k(t_1)  v_k + \sum_{j = 1}^N \sum_{k = 1}^j (\psi_k(t_1) \Delta\psi_j(t_1) + \psi_j(t_1) \Delta\psi_k(t_1))h^{jk} + \dots.
\end{equation}
Under the flow of Equation \eqref{maineq}, $\psi_k(t_2) = \psi_k(t_1) \exp(\lambda_k (t_2-t_1))$ and $\Delta \psi_k(t_2) = \Delta \psi_k(t_1) \exp(\lambda_k (t_2-t_1))$.  As such, at any other time $t_2$  one can write
\begin{align}
\Delta x(t_2) &= \sum_{k = 1}^N \Delta \psi_k(t_1) \exp(\lambda_k (t_2-t_1))  v_k   \nonumber \\
& \quad +\sum_{j = 1}^N \sum_{k = 1}^j (\psi_k(t_1) \Delta\psi_j(t_1) + \psi_j(t_1) \Delta\psi_k(t_1)) \exp((\lambda_k + \lambda_j)(t_2-t_1))h^{jk} + \dots.
\end{align}
Assuming that $|\lambda_k| + |\lambda_j| <  |\lambda_\beta|$ for all $k,j \leq \beta$, in the limit that that $t$ becomes large, one finds
\begin{equation} \label{xteq}
    \Delta x(t_2) = \sum_{k = 1}^\beta \Delta \psi_k(t_1) \exp(\lambda_k (t_2-t_1))  v_k + O(\epsilon),
\end{equation}
where $0 < \epsilon \ll \exp(\lambda_k (t_2-t_1)) \leq 1$ for $k \leq \beta$ (recall that ${\rm Real}(\lambda_k) \leq 0$ for all $k$).  $\Delta x(t_2)$ can also be computed as the solution to the linear time varying equation $\Delta \dot{x} = J(t) \Delta x$ where $J$ is the Jacobian evaluated at $x(t)$.  The solution of linear time varying system is 
\begin{equation} \label{statetransition}
    \Delta x(t_2) = \Phi_J(t_2,t_1) \Delta x(t_1),
\end{equation}
where $\Phi_J(t_2,t_1)$ is the state transition matrix.  For clarity of exposition, assume $\Phi_J(t_2,t_1)$ is diagonalizable.  Comparing Equation \eqref{statetransition} to \eqref{xteq} the order 1 terms of $\Delta x(t_2)$ can be written as a linear vector space spanned by $\{v_1,\dots,v_\beta\}$.  With this in mind, letting $\hat{\lambda}_k$ be an eigenvalue of $\Phi_J(t_2,t_1)$ with corresponding left and right eigenvalues $\hat{w}_k$ and $\hat{v}_k$, respectively.  After diagonalization, $\Phi_J(t_2,t_1)$ can be well approximated by
\begin{equation}  \label{phijeq}
\Phi_J(t_2,t_1)  \approx  \underbrace{\sum_{k = 1}^\beta  \hat{v}_{k}  \hat{\lambda}_{k} \hat{w}_{k}^T}_{O(1)\; {\rm terms}}  +  \underbrace{ \sum_{k = \beta  +1}^N  \hat{v}_{k}  \hat{\lambda}_{k} \hat{w}_{k}^T}_{O(\epsilon)\; \rm{terms}},
\end{equation}
 where $\{\hat{v}_1,\dots,\hat{v}_\beta\}$ form a basis for the span of $v_1,\dots,v_\beta$.  Additionally, because the state transition matrix is always invertible $v_k \notin {\rm Span}(\{\hat{v}_1,\dots,\hat{v}_\beta\})$ for $k>\beta$. Note that a similar result can be obtained in the case where $\Phi_J(t_2,t_1)$ is not diagonalizable can also be handled by considering Jordan blocks of size larger than 1.

\section{Asymptotic Expansion of the State in a Basis of Isostable Coordinates} \label{apxb}
\renewcommand{\thetable}{B\arabic{table}}  
\renewcommand{\thefigure}{B\arabic{figure}} 
\renewcommand{\theequation}{B\arabic{equation}} 
\setcounter{equation}{0}
\setcounter{figure}{0}

In \cite{wils21dd}, a strategy was described for computing the terms $h^{ij}$, $h^{kjk}$, $\dots$ in the expansion \eqref{xexp} for $x$ in a basis of isostable coordinates.  A description of this approach is provided here for convenience.  

To begin, letting $F(x) = \begin{bmatrix} f_1 & \dots & f_N \end{bmatrix}^T$, defining $\Delta x = x-x_0$ and asymptotically expanding \eqref{maineq} about the fixed point yields
\begin{align} \label{taylfn}
\frac{d \Delta x}{dt} &= J   \Delta x  +  \begin{bmatrix}  \sum_{i = 2}^\infty \frac{1}{i !} \left[ \overset{i}{\otimes}  \Delta x^T   \right]{\rm vec}( f_1^{(i)})  \\ \vdots \\ \sum_{i = 2}^\infty \frac{1}{i !} \left[ \overset{i}{\otimes}  \Delta x^T   \right] {\rm vec}( f_N^{(i)})     \end{bmatrix},
\end{align}
where $\otimes$ is the Kronecker product, ${\rm vec}(\cdot)$ is an operator that stacks each column of a matrix into a single column vector, and higher order partial derivatives are defined recursively as
\begin{equation} \label{fkeq}
f_j^{(k)} = \frac{\partial {\rm vec} \big( f^{(k-1)}_j \big)}{\partial x^T} { \in \mathbb{R}^{ N^{(k-1)} \times N}},
\end{equation}
where all partial derivatives are evaluated at $x_0$.  The notation $ \overset{3}{\otimes}  \Delta x^T$  is shorthand for $\Delta x^T \otimes \Delta x^T \otimes \Delta x^T {\in \mathbb{R}^{1 \times N^3}}$.  Equating \eqref{deriveq} and \eqref{taylfn}, substituting $\Delta x$ from \eqref{xexp} into \eqref{taylfn}, and matching the resulting powers in the isostable coefficients yields the following relationship
\begin{equation} \label{geqs}
0 = (J - ( \lambda_i + \lambda_j + \lambda_k +\dots ) {\rm Id} ) h^{ijk\dots} + q^{ijk\dots},
\end{equation} 
where ${\rm Id}$ is an appropriately sized identity matrix and $q^{ijk\dots} \in \mathbb{R}^N$ is a function of the lower order terms of the expansion. For example, $q^{321}$ (a third order term) might depend on $v_1$ or $h^{21}$ (a first or second order term) but will not depend on $h^{111}$ (a third order term).  As such, $g^{ijk\dots}$ can be computed as the solution to a linear matrix equation starting with lower order terms and progressing to higher order ones.  From a computational perspective, note that the size of the matrix $f_j^{(k)}$ in Equation \eqref{fkeq} grows exponentially with increasing $k$.  In previous works using similar methods \cite{wils21input}, \cite{wils20highacc}, computation of \eqref{xexp} beyond $8^{\rm th}$ order accuracy becomes computationally infeasible when working with models of dimension $N>2$. Additionally, diminishing returns are typically seen when \eqref{xexp} is taken to increasingly large orders of accuracy.

\end{appendices}

\end{document}